
\input amstex
\def\aa#1{A_{#1}}
\def\aaa#1#2{A_{#1}^{#2}}
\def\al#1#2#3{a^{#1}_{{#2}{#3}}}

\def\bii#1#2{\bold F\bold F^{#1}_{#2}}
\def\bcc#1{\Bbb C^{#1}}
\def\bpp#1{\Bbb P^{#1}}
\def\bppp{\Bbb P }\def\bccc{\Bbb C }\def\baaa{\Bbb A }
\def\bb#1#2#3{b^{#1}_{{#2}{#3}}}

\def\cf{\Cal F}
\def\ci{\Cal I}

\def\ff#1#2{F_{#1}^{#2}}

\def\gg#1#2{g^{#1}_{#2}}

\def\hd{, \hdots ,}
\def\inv{{}^{-1}}
\def\ii{| II_X|}

\def\ooo#1#2{\omega^{#1}_{#2}}
\def\oo#1{\omega^{#1}}
\def\ot{\!\otimes\!}

\def\pp#1{\Bbb P^{#1}}
\def\ppp{\Bbb P}
\def\qq#1#2#3{q^{#1}_{{#2} {#3}}}
\def\qqq#1#2#3{q^{n+ #1}_{{#2} {#3}}}
\def\ra{\rightarrow}
\def\rr#1#2#3#4{r^{#1}_{{#2} {#3}{#4}}}
\def\schur{S^{(21)}_{\circ}}
\def\tim{\text{Image}\,}
\def\tdim{\text{dim}\,}
\def\tker{\text{ker}\,}
\def\tmod{\text{ mod }}

\def\trank{\text{rank}\,}
\def\ud#1{\underline d^{#1}}
\def\up#1{{}^{({#1})}}
\def\upperp{{}^{\perp}}
\def\vdx{v_d(X)}

\def\vtwox{v_2(X)}
\def\vthreex{v_3(X)}

\def\ww{\wedge}
\def\xx#1{x^{#1}}
\def\xsm{X_{sm}}

\define\intprod{\mathbin{\hbox{\vrule height .5pt width 3.5pt depth 0pt %
        \vrule height 6pt width .5pt depth 0pt}}}
\documentstyle{amsppt}
\magnification = 1200
\hsize =15truecm
\hcorrection{.5truein}
\baselineskip =18truept
\vsize =22truecm
\NoBlackBoxes
\topmatter
\title
Differential-Geometric Characterizations of Complete Intersections \endtitle
\rightheadtext{Differential-Geometric Characterizations}
\author
  J.M. Landsberg
\endauthor

\address{Department of Mathematics,
Columbia University, New York,  NY 10027}
\endaddress
\email {jml\@math.columbia.edu }
\endemail
\date {  June 25, 1995 }\enddate
\thanks {This work was done while the author was  supported by an NSF
postdoctoral fellowship and later by NSF grant DMS-9303704.}
\endthanks
\endkeywords
\subjclass{ primary 14e335, secondary 533a20}\endsubjclass
\abstract{ We  characterize complete intersections
in terms of local differential geometry.

Let $X^n\subset\Bbb C\Bbb P^{n+a}$ be a variety. We first
localize the problem; we give a criterion for $X$ to
be a complete intersection that is testable at any smooth point of
$X$. We rephrase the criterion in the language of projective differential
geometry and derive a sufficient condition for $X$ to be a complete
intersection that is computable at a general point $x\in X$.
The sufficient condition has a geometric interpretation in terms
of restrictions on the spaces of osculating hypersurfaces at $x$.
 When this sufficient
condition holds, we are able to define
systems of partial differential equations
that generalize  the classical
Monge equation that characterizes conic curves in $\bccc\Bbb P^2$.

  Using our sufficent condition, we show that if
the ideal of $X$ is generated by quadrics  and
 $a<\frac{n-(b+1)+3}{3}$, where $b=$dim$X_{sing}$,
 then $X$ is a complete intersection.
 }
\endabstract

\endtopmatter

\document

\noindent\S 0. {\bf Introduction}

\smallpagebreak

\noindent{\bf Local and global geometry}

Projective differential geometry has been used to study
the local geometry of subvarieties of projective space
by various authors
(e.g. [C], [F], [GH], [JM], [T]).  However, there are few examples where
global conclusions are drawn from the local picture.

One (global) fact  about projective varieties
that has been encoded into the infinitesimal geometry is the following: if
there is a line on a
variety $X^n\subset\bccc\bpp{n+a}$ along
 which the embedded tangent space is constant, then $X$ must
be singular.  Griffiths and Harris realized this
fact  had implications
for the projective second fundamental form (that its singular locus
must be empty at general points) which  enabled them
to
reprove Severi's theorem that the only smooth
surface in $\bpp 5$ with degenerate secant variety is the Veronese ([GH],
6.18),
using local methods.

In [L1], we used Zak's theorem on tangencies ([Z], 1.8), or equivalently the
Fulton-Hansen Connectedness theorem
([FL], 3.1) to encode additional global information
about projective varieties into the local
differential geometry. We proved a {\it rank restriction theorem}
([L1], 6.1) that bounds from below the ranks of quadrics in  the
projective second fundamental form. A general  principle
that the rank restriction theorem illustrates is:

\smallpagebreak

{\it
In order for a variety of small codimension to be smooth, it
must \lq\lq bend enough\rq\rq .}

\smallpagebreak

In the case of the rank restriction theorem, \lq\lq bending
enough\rq\rq\  corresponds to genericity of the projective second
fundamental form at general points.

The rank restriction theorem
 enabled us to reprove Zak's theorem on Severi varieties
and to prove further results about varieties with degenerate secant
varieties using local methods  in [L2].
One consequence of the rank restriction
theorem that is particularly easy to apply is:

\proclaim{Theorem [L1], 4.14} Let $X^n\subset\bccc\bpp{n+a}$ be a  variety with
$a<\frac
{n-(b+1)}2+1$
(where $b=\tdim X_{sing}$).
Then at general points of $X$, the third fundamental form of $X$ is zero.
  \endproclaim

In this paper we apply the results of [L1] to obtain information about the
quadrics containing a
variety.  The   results on quadrics  follow from an in-depth local study that
occupies
the bulk of this paper.

The local study shows precisely how one can detect  the failure
of a variety to be a complete intersection from the differential invariants
at a general point $x\in X$.
Vaguely stated, the motivating
principle  is:

\smallpagebreak

{\it  If $X$ is
not a complete intersection, it \lq\lq bends less\rq\rq\ than
expected.}

\smallpagebreak

 Work of others (e.g. [R]) on complete intersections
can be understood in these terms
when one uses a K\" ahler metric. In this paper, \lq\lq bending less\rq\rq\
will correspond to
certain  non-generic behavior of projective differential
invariants and the \lq\lq expectation\rq\rq\
 from knowledge about the generators
of the ideal of $X$ that is also computable from the projective differential
invariants.

\medpagebreak

\noindent{\bf Osculating hypersurfaces}

We use  differential invariants to study the hypersurfaces containing $X$ via
the {\it osculating hypersurfaces} at a  point $x\in X$. Osculating
hypersurfaces
 are geometric objects
of interest in their own right.
 There are several equivalent definitions of what it means for a hypersurface
to
osculate to order $k$ at a point $x\in X$,  an elementary one is as follows:

\smallpagebreak

\noindent{\bf Definition 0.1} Let $X^n\subset\bpp{n+a}$ be a variety and let
$Z\subset\bpp{n+a}$ be a hypersurface. Let $x\in X\cap Z$.
  In an affine open set $\baaa$ containing $x$, $Z$
is given locally by  a function $f$ which restricts  to a function $\overline
f$
on $X\cap \baaa$. If we identify the completion of the local ring of $X\cap
\baaa$ at $x$ with
$\bcc{}[[ \xx 1\hd \xx n]]$,
 then $Z$ {\it osculates to order  $k$ at $x$}  if
the power series of $\overline f$  at $x$ has no terms of degree
less than or equal to $k$.

\smallpagebreak

 Let $V=\bcc{n+a+1}$ and let $X^n\subset\bppp V=\bcc{}\bpp{n+a}$ be a variety
of dimension $n$. Let $I_X\subset S^{\bullet}V^*$ denote
the ideal of $X$
and let $I_{X,d}=I_d = S^dV^*\cap I_X$ denote
the $d$-th graded piece of $I_X$.

$I_d$,
the vector space of all hypersurfaces of degree
$d$
containing $X$,
can be characterized geometrically as
follows: Let
$\vdx\subset\bppp S^dV$ denote the $d$-th Veronese re-embedding of $X$,
and let $<\vdx >$ denote its linear span.
 $ I_d ={<\vdx >\!\upperp\! }
\subset  S^dV^*$,    the
annhilator of
the linear span of $\vdx$.
Similarly, the space of hypersurfaces of degree $d$ osculating to order $k$
at a smooth point $x\in X$ is
$ \tker\bii k{\vdx x}\subset S^dV^*$,    the kernel of the
{\it  $k$-th fundamental form} of $\vdx$ at $x$. (See \S 2 for definitions of
the fundamental forms.)

\smallpagebreak

Using the language of osculating hypersurfaces, [L1],(4.14)   may be restated
as:

\proclaim{Restatement of [L1], 4.14}Let $X^n\subset\bccc\bpp{n+a}$ be a
variety with $a<\frac
{n-(b+1)}2+1$ (where $b=\tdim X_{sing}$). Let $x\in X$ be a general point.
 If a hyperplane  $H$ osculates to
 order two at  $x$, then $X\subset H$.
\endproclaim

In contrast,
 any  surface in  $\bpp 6$  has at least one hyperplane osculating to order two
at  {\it every } point.  There is also a class of   smooth
surfaces in $\bpp 5$, the {\it Legendrian
surfaces}, which have the property that
at every point there is at least one hyperplane osculating to order two. (This
class   includes the ruled
surfaces.)

\medpagebreak

\noindent{\bf Overview}

 The precise meanings of terms used
in this overview are explained in the sections.

\smallpagebreak

In \S 1 we observe that if $X$ is such that $I_X = (I_d)$ and
$I_{d-1}=(0)$, then $X$ is a complete intersection if and only if
every  hypersurface $Z\in I_d$ is smooth at all
the smooth points of $X$ (1.1). We generalize (1.1) to a local characterization
of all complete intersections.  The characterization is simplified
thanks to the following definition due to L'vovsky [L'v]:

\noindent{\bf Definition 1.2}. Let $X\subset\ppp V$ be a variety. Let
$P\in I_d$ and let $Z=Z_P\subset\ppp V$ be the corresponding hypersurface.
We will say
$Z$ {\it trivially contains} $X$ if
$P= l^1P_1+\hdots +l^kP_k$ with $P_1\hd P_k\in I_{d-1}$ and $l^1\hd l^k
\in V^*$, and   otherwise
that  $Z$ {\it essentially contains} $X$.

 \proclaim{Proposition 1.6,
A local characterization of complete intersections} Let $X\subset\bppp V$ be a
variety and let $\xsm$ denote its smooth points. The following are equivalent:

1. $X$ is a complete intersection.

2. Every hypersurface
essentially containing $X$ is smooth at all $x\in\xsm$.

3.  Let $x\in \xsm$. Every
 hypersurface
essentially containing $X$ is smooth at $x$.

4.  Let $x\in \xsm$. For all $k$, the map
$$
\align
[d]_k : I_k/(I_{k-1}\circ V^*)  &\ra  N^*_xX/ d(I_{k-1}\circ V^*)\\
[P] &\mapsto [dP_x]
\endalign
$$
is injective. ($dP$ denotes the exterior derivative of the
polynomial $P$.)

 \endproclaim

In fact, (1.6) not only localizes the question of whether or not
$X$ is a complete intersection, it allows one to study the property
one degree at a time.

\smallpagebreak

\noindent{\bf Definition 1.8}. Fix a point  $x\in \xsm$.
We will say
{\it $X$ has no excess equations in degree $k$ at $x$},
or  that  {\it $(CI)_k$ holds at $x$}, if $[d]_k$ is injective
at $x$. (Note that $X$ is a complete intersection if and only if
$(CI)_k$ holds at $x$ for all $k$).

\smallpagebreak

Unfortunately,  to determine if $X$ satisfies
$(CI)_k$ at a point $x$ is not necessarily computable in a predictable
number of steps. To avoid this problem, we restrict attention to
sufficient conditions for a variety to be a complete intersection
in the remainder of the paper. The sufficient conditions are expressed
geometrically in terms of the osculating hypersurfaces of $X$ at
$x$, and to deal with osculating hypersurfaces we need to deal
with the projective differential invariants of $X$.

\smallpagebreak

In \S 2, we review the {\it projective fundamental forms}
of a variety,
denoted $\bii k X$ in this paper.
 Of special importance is $\bii 2 X=II_X$, the
projective second fundamental form. We also need to work with
some subtler differential invariants, which we denote
$F_k$ and call the {\it $(k-2)$-nd variation of $II$}.

\smallpagebreak

In \S 3, we study the
space of  hypersurfaces of
degree $d$ osculating to order $k$   at a smooth point $x\in X$,
taking advantage of their description as $\tker\bii k{\vdx ,x}$ mentioned
above. We compute the fundamental forms of the Veronese re-embeddings
of $X$ and prove some generalities about the osculating
spaces.
 First we show that the dimension of the space of hypersurfaces of
degree $d$ osculating to order $k$ at $x$ is fixed for all
$k\leq d$, generalizing the fact that there is an
$(a-1)$-dimensional space of hyperplanes osculating to order
one at each smooth point:

 \proclaim{Proposition  3.16} Let $X^n\subseteq\bpp{} V=\Bbb C\pp{n+a}$ be a
variety and let  $x\in \xsm$.
 For all $  p\leq d$,
$$
\align &
\tdim
\left\{
\matrix  \text{(not necessarily irreducible)
hypersurfaces} \\\text{of degree }d\text{ osculating   to order }  p  \text{
at }x\endmatrix \right\} \\
&
= \binom{(n+a+1)+(d-1)}  d -  \{1+
n+ \binom{ n+1} 2  + \hdots +
\binom{ n+p-1}p
 \}.\endalign
$$
\endproclaim

 For $k>d$, the dimensions depend
on the geometry of $X$. For $d+1\leq k\leq 2d-1$ there are {\it lower} bounds
on the dimensions
of the space of hypersurfaces of
degree $d$ osculating to order $k$ at $x$. For example:

 \proclaim{Proposition  3.17} Let $X^n\subseteq\bpp{} V=\Bbb C\pp{n+a}$ be any
variety,
 and $x\in X $, any smooth point.
$$
\tdim
\left\{
\matrix  \text{(not necessarily irreducible )
hypersurfaces of}\\\text{degree }d\text{  osculating   to order } 2d-1  \text{
at }x\endmatrix \right\}
  \geq
\binom {a+d-1}d-1.
$$
 \endproclaim

We rephrase (1.6) in a manner suitable for
approximation by osculating hypersurfaces.

We include several examples, including a review of the
classical Monge equation; a fifth order differential equation
that characterizes conic curves in $\pp 2$.

 We define a
condition which we denote
$(CI)^{2k}_k$,
that implies $(CI)_k$ holds (3.22). If the codimension of $X$ is
sufficiently small,  $(CI)^{2k}_k$ is a genericity condition.
While only being a sufficient condition for $(CI)_k$,
$(CI)^{2k}_k$
has the advantage of being computable by taking
at most $2k$ derivatives.

 Finally we
describe generalizations of the classical Monge equation. The classical
Monge
equation actually only classifies
\lq\lq nondegenerate\rq\rq\ conic curves, those that are not
pairs of lines. The generalized Monge systems described by (3.23)
characterize \lq\lq nondegenerate\rq\rq\
complete intersections  whose ideals are generated
in degrees $d_1<\hdots <d_r$, by a PDE system of order
$2d_r+1$. In this case \lq\lq nondegenerate\rq\rq\
means that the   conditions  $(CI)^{2{d_1}}_{d_1}\hd
(CI)^{2{d_r}}_{d_r}$ hold, which is a
natural
genericity condition in small
codimension.

\smallpagebreak

In \S 4 we study the osculating
quadrics of $X$ in detail.  We utilize
a condition slightly stronger than
$(CI)^4_2$,
 which we call   {\it strong genericity
in degree two}. Strong genericity
in degree two is
computable by taking two derivatives at a general point. It is
  a  genericity
condition on $II_{X,x}$; namely that the system of quadrics
$|II_X|_x$  has
no linear syzygies. Our results are as follows:

We define a  more precise generalized Monge system for complete
intersections of quadrics (4.17), in the sense that (4.17) is expressed
in terms of the invariants $F_k$ instead of
 the fundamental
forms of $\vtwox$.

Theorem (4.18)  states that if there are no linear syzygies in
$II_X$ at general points, then there are {\it upper} bounds on the dimensions
of the spaces of quadrics osculating to orders three and four,
complementing (3.17).  Moreover, these
bounds are achieved if and only if the
(fifth order)  generalized Monge system  for quadrics holds, and the
only way such varieties can be cut out by quadrics is if the generalized
Monge system holds.

A consequence of (4.18) is:

 \proclaim{Corollary 4.20} Let $X\subset\bppp V$ be a variety and $x\in X$ a
general point.
Assume $III_{Xx}=0$. If   there are no quadric hypersurfaces singular at $x$
that osculate to order four at $x$, and $Q$ is a quadric hypersurface
osculating to order five at $x$, then $X\subseteq Q$.\endproclaim

Theorem (4.28) describes a sufficient condition for $I_X$ to be
generated by quadrics.

Theorem (4.33) describes the structure of varieties whose
projective differential invariants of order
greater than two are zero at general points.

We observe that if $III_X=0$ and if the ideal of
$X$ is generated by quadrics, it is
extremely difficult for $X$ to fail to be a complete intersection without being
an intersection of a complete intersection
with a  rational variety (4.39).  We illustrate this principle  with some
examples.

\smallpagebreak

In \S 5 we describe the equations of some familiar homogeneous varieties in
a manner that illustrates the results in \S 4. Among the
varieties treated are
 the Severi varieties and the ten-dimensional spinor
variety, $\Bbb S^{10}\subset\pp{15}$.

\smallpagebreak

In \S 6 we study systems of quadrics with linear syzygies and systems with
nonzero prolongations.
In representation-theoretic language the study is as follows:
Let  $A\subset S^2T^*$  be a system of
quadrics.
We determine rank restrictions  that force the intersection of $A\ot T^*$ with
the two
irreducible
$Gl(T^*)$ modules in $S^2T^*\ot T^*$ to be the origin.
We  combine the results here with the results of [L1] and
\S 4 to prove:

\proclaim{Theorem 6.26 }
Let  $X^n\subset\bpp{n+a}$ be a variety and $x\in X$ a general point.
  Let $b=\tdim X_{sing}$. (Set $b=-1$ if $X$ is smooth.)
If
$a<\frac{n-(b+1)+3}{3}
$
then
$$
\align &
\tdim \{ \text {quadrics osculating to order three at } x \}\leq
a
+\binom{a+1}2-1 \tag 6.27\\
&
\tdim \{ \text {quadrics osculating to order four at } x \}\leq
a-1.
\endalign
$$
 Equality occurs in the first (respectively second)  line of  (6.27)
 if and only if
(4.9.3) (resp. (4.9.4))   holds at $x$.
  If  the generalized Monge system
(4.17)  holds, then  $X$ is a complete intersection of the
$(a-1)$-dimensional family of quadrics osculating to order four.
\endproclaim

\proclaim{Corollary 6.28}Let  $X^n\subset\bpp{n+a}$ be a variety and
 $x\in X$ a general point.
  Let $b=\tdim X_{sing}$. (Set $b=-1$ if $X$ is smooth.)
If
$a<\frac{n-(b+1)+3}{3}
$
then any quadric osculating to order four at $x$ is smooth at $x$ and
any quadric osculating to order five at $x$ contains $X$.
\endproclaim

\proclaim{Corollary 6.29 }
Let  $X^n\subset\bpp{n+a}$ be a variety
with $I_X$ generated by quadrics.   Let $b=\tdim
X_{sing}$. (Set $b=-1$ if $X$ is smooth.)
If
$a<\frac{n-(b+1)+3}{3},
$
then $X$ is a complete intersection.
\endproclaim

\smallpagebreak

\noindent{\bf Acknowledgements.}
The inspiration for this paper came from a series of
conversations with Mark Green, who has also provided
 substantial help at each
step along the way. After reading a preliminary version of this paper,  S.
L'vovsky pointed out an error in Lemma 6.19 and how to correct it. He also
came up with some new terminology which is vastly superior to the
terminology that was in the old version.   The referee has
also provided numerous suggestions and corrections which have hopefully
made the paper more readable. It is a pleasure to thank
professors Green, L'vovsky, and the referee for their help.

\smallpagebreak

\noindent{\bf Notation}. We  will use the following conventions for indices
$$
\align &0\leq B,C\leq n+a \\
& 1\leq \alpha ,\beta \leq n \\
& n+1\leq \mu ,\nu \leq n+a   \endalign
 $$

Alternating products of  vectors will be denoted with a wedge ($\wedge$), and
symmetric
products will not have any symbol (e.g. $\omega\circ\beta$ will be denoted
 $\omega\beta$).
$T_xX$ denotes the holomorphic tangent space to $X$ at $x$
and $\tilde T_xX$ the embedded tangent space. In
general we will supress reference to the base point of our manifold $X$ when we
abbreviate the names of bundles
so $T$ should be read as $T_pX$ for a general $p\in X$, $N$ as $N_pX$ etc... .
  $\{ e_i \}$ means the span of the vectors $e_i$ over the index range $i$.
If $Y\subset\bpp m$ then $\hat Y\subset\Bbb C^{m+1}$ will
be used to  denote the cone over
$Y$
 (with the exception that  the cone over the embedded tangent space
${\tilde T}$ will be denoted
$\hat T$).
We will often ignore twists in bundles so $T$ will
be used to denote both
$T_xX$ and $T_xX(1)=\hat T/\hat x$.
If $A\in \Bbb C^{n+a+1}$, its projection to $\bpp{n+a}$ will be denoted $[A]$.
If $V$ is a vector space and $W$ a subspace, and $(e_1\hd e_n )$ a basis of $V$
such that $\{e_1\hd e_p \}=W$, we write
$\{e_{p+1}\hd e_n\}$ mod$W$ to denote the space $V/W$.
For vector subspaces $W\subset V$, we will use the notation
$W\upperp  \subset V^*$ for the annihilator of $W$ in $V^*$.
We will use the summation convention throughout (i.e. repeated indices are to
be summed over).
$\frak S_{\alpha\beta\gamma}$ denotes cylic summation over
the fixed indices $\alpha\beta\gamma$.
In general, $X$ will denote a variety,  $\xsm$ its smooth points, and
$X_{sing}$ its singular points.  $\Bbb C\bpp k$ will  be denoted $\bpp k$.
$\bii k X$ is the $k$-th fundamental form of $X$.
We will often denote $\bii 2X$ by $II$ and $\bii 3X$ by $III$.
$F_k=F_k^{\mu}A_{\mu}\tmod \hat T$ is the differential invariant we
call the  {\it $(k-2)$-nd variation of $II$}.
By a {\it general point}
$x\in X$ we mean  a smooth point of $X$ such that all the discrete information
in the
differential invariants of $X$ is locally constant.  The nongeneral points of
$X$ are
a codimension one subset of $X$.

 \bigpagebreak

\noindent\S 1. {\bf An elementary characterization of complete intersections}

 Let $V=\bcc{n+a+1}$ and let $X^n\subset\bppp V=\bcc{}\bpp{n+a}$ be a variety
of dimension $n$.
Let $\xsm$ denote the smooth points of $X$.
 Let $I_X\subset S^{\bullet}V^*$ denote
the ideal of $X$
and let $I_{X,d}=I_d = S^dV^*\cap I_X$ denote
the $d$-th graded piece of $I_X$. Fixing a smooth point  $x\in X$,
there is a distinguished subspace of  $I_d$, namely the
hypersurfaces of degree $d$ that are singular at $x$, i.e. $P\in I_d$ such that
$(dP)_x=0$, where $dP$ denotes the
exterior derivative of the polynomial  $P$.

\proclaim{Proposition 1.1} Let $X\subset\bppp V$ be a variety such
that  $I_X= (I_d)$ (i.e. $I_X$ is generated by
$I_d$)
and $I_{d-1}=(0)$. Then the following are equivalent:

1.  $X$ is a complete intersection.

2. Every hypersurface of degree $d$ containing $X$ is smooth at all
$x\in\xsm$.

3.  Let $x\in \xsm$. Every
 hypersurface of degree $d$ containing $X$ is smooth at $x$.

 \endproclaim
\demo{Proof}
Say $P\in I_d$ is nonzero and such that $dP_x=0$. Let $x
\in \xsm$ be a smooth point, then
there
exist  $P_1\hd P_a\in I_d$ such that
$(dP_1)_x\hd (dP_a)_x$ span
the conormal space $N^*_xX\subset T^*_x\bppp V$
(actually $N^*_xX(1)\subset T^*_x\bppp V(1)$).
Since
$(dP_1)_x\hd (dP_a)_x$ are linearly independent, $P$ is not in the ideal
generated
by$ P_1\hd P_a $, and
thus $X$ is not a complete intersection.
Conversely, if $X$ is not a complete intersection, one can always find such a
$P$.
\qed\enddemo

The following definition is due to L'vovsky [L'v]:

\noindent{\bf Definition 1.2}. Let $X\subset\ppp V$ be a variety. Let
$P\in I_d$ and let $Z=Z_P\subset\ppp V$ be the corresponding hypersurface.
We will say
$Z$ {\it trivially contains} $X$ if
$P= l^1P_1+\hdots l^mP_m$ with $P_1\hd P_m\in I_{d-1}$ and $l^1\hd l^m
\in V^*$, and   otherwise
that  $Z$ {\it essentially contains} $X$.

\smallpagebreak

Note that the space of hypersurfaces
of degree $k$ that trivially contain
$X$   is $ I_{k-1}\circ V^*$. Fix $x\in \xsm$ and consider the map
$$
\align
[d]_k : I_k/(I_{k-1}\circ V^*)  &\ra  N^*_xX/ d(I_{k-1}\circ V^*)\tag 1.3.k\\
[P] &\mapsto [dP_x]
\endalign
$$
where $dP$ denotes the exterior derivative of the polynomial $P$.
(Again, here and in what follows, we really should be writing
$N^*_xX(1)$.)

\smallpagebreak

\noindent{\bf Definition 1.4}. Fix $x\in \xsm$.
Let $N_k^*:=\{ dP_x | P\in I_k \}\subseteq N^*_xX$.
Let $d_1\hd d_r$ be the smallest integers such that
$$
0\subsetneq N_{d_1}^*\subsetneq N_{d_2}^*\subsetneq
\hdots\subsetneq N_{d_r}^*=N^*_xX.\tag 1.5
$$
 We will
call   (1.5)
 {\it the natural filtration of $N^*_xX$}.

\smallpagebreak

(1.1) generalizes to the following
statement:

 \proclaim{Proposition 1.6,
A local characterization of complete intersections} Let $X\subset\bppp V$ be a
variety. The following are equivalent:

1. $X$ is a complete intersection.

2. Every hypersurface
essentially containing $X$ is smooth at all $x\in\xsm$.

3.  Let $x\in \xsm$. Every
 hypersurface
essentially containing $X$ is smooth at $x$.

4.  Let $x\in \xsm$   be any smooth point of $X$. For all $k$, the map
$$
\align
[d]_k : I_k/(I_{k-1}\circ V^*)  &\ra  N^*_xX/ d(I_{k-1}\circ V^*)\\
[P] &\mapsto [dP_x]
\endalign
$$
is injective. ($dP$ denotes the exterior derivative of the
polynomial $P$.)

 \endproclaim

\demo{Proof}

Consider the sum of all the maps (1.3.k),
$$
[d]^{\oplus} : \oplus_k I_k/(I_{k-1}\circ V^* ) \ra \oplus_k (N^*_k/N^*_{k-1})
\tag 1.7
$$
where $N_k^*$ is as in (1.4)

The dimension of the target of $[d]^{\oplus}$ is $a$,
and $[d]^{\oplus}$ is surjective so  the dimension
of the source is $a$ if and only if $[d]^{\oplus}$ is injective, i.e.
all the maps (1.3.k) are injective. On the other hand, the
dimension of the source is exactly the number of polynomials
needed for a minimal set of generators of $I_X$, and the
kernel of $[d]^{\oplus}$ is exactly the hypersurfaces that
essentially contain $X$ and are singular at $x$.
\qed\enddemo

\noindent{\bf Definition 1.8}. Fix a point  $x\in \xsm$.
We will say
{\it $X$ has no excess equations in degree $k$ at $x$},
or  that  {\it $(CI)_k$ holds at $x$}, if (1.3.k) is injective.
The equivalence $1\Leftrightarrow 3$ of (1.6) may
be rephrased as:
$X$ is a complete intersection iff $(CI)_k$ holds  for all $k$
at some $x\in \xsm$.

Although the following is clear from the discussion above,
we record it as a proposition.

\proclaim{Proposition 1.9 } Let $X\subset\bppp V$ be a variety and
$x\in X$ a smooth point.
  Let  $N_{d_1}^*\subset N_{d_2}^*\subset\hdots\subset N_{d_r}^*=N^*$ be
 the natural filtration of $N^*_xX$ described in (1.4) and let
$a_j= \tdim (N^*_{d_j} /N^*_{d_{j-1}})$. Then there  is
at least an
$(a_1-1)$-dimensional space of irreducible hypersurfaces of degree $d_1$
essentially containing $X$,  an
$(a_2-1)$-dimensional space of  irreducible hypersurfaces of degree $d_2$
essentially containing $X$, ...., and an
$(a_r-1)$-dimensional space of irreducible  hypersurfaces of degree $d_r$
essentially containing $X$.

 $X$ is a complete intersection   if and only if $I_X$ is   generated by the
equations for these hypersurfaces.
 \endproclaim

So far we have localized the study of complete intersections to
a point, and further, filtered the conormal bundle at that point
to enable us to study one degree at a time. Unfortunately, to determine
if a hypersurface essentially contains $X$, one might need to take
an arbitrarily high number of derivatives. To have computable
conditions, we will work  with osculating
hypersurfaces rather than the
hypersurfaces
 containing $X$. The advantage will be that we will only
need to study a fixed number of derivatives for each fixed degree of
hypersurface;  the disadvantage is that we will only obtain
sufficient conditions to be a complete intersection.

 We first review some
notions from projective differential geometry.

 \bigpagebreak

\noindent\S 2. {\bf Frames, fundamental forms}

\smallpagebreak

Refer to   [GH],[L1] for  more detailed explanations of what follows.

Let $X\subset\bpp{} V =\bcc{}\bpp{n+a}$ be a variety and $x\in X$ a general
point.
We will  compute fundamental forms using   a method described in [GH]
which will be  useful  for the computations of this paper.
One may take what follows as the definitions of the fundamental forms,
although more
geometric definitions are given in [L1].

Given $X\subset\bpp{} V$ and $x\in \xsm$, let $\hat x\subset V$ denote the line
$x$ determines and let $\aa 0 \in \hat x$.
Extend $\aa 0$ to a basis  $f
=(\aa 0\hd\aa{n+a})=(\aa 0,\aa\alpha ,\aa\mu )$ of $V$ adapted
to the filtration $\hat x\subset\hat T\subset V$, where
$1\leq\alpha ,\beta\leq n,\  n+1\leq\mu,\nu\leq n+a$, so
$\hat T
=\{\aa 0,\aa\alpha \}$ is the cone over the embedded tangent space to $X$ at
$x$.
  Let
$\cf^1\ra \xsm$ be the bundle of all such bases.

The {\it Maurer-Cartan form}  of $Gl(V)$, defined by
 $$
\Omega_f:=
f\inv df \tag 2.1
$$
restricts to $\cf^1\subset Gl(V)$ and satisfies the
{\it Maurer-Cartan equation}
$$
d\Omega = -\Omega\ww\Omega .\tag 2.2
$$
See ([GH],[L1]) for more on the Maurer-Cartan form.

 $\cf^1$ is a principal $G_1$-bundle over $\xsm$ where
$$
G_1 := \left\{ g\in Gl(V) \mid g= \pmatrix g^0_0 & g^0_{\beta} & g^0_{\nu} \\
0&g^{\alpha}_{\beta}&g^{\alpha}_{\nu} \\  0&0&g^{\mu}_{\nu}\endpmatrix\right\}.
 \tag 2.3
$$

Define inductively a series of maps:
$$
\ud k\aa 0 : (T\cf^1)^{\ot k}\ra V/\tim
(\ud 0\hd \ud{k-1}) \tag 2.4
$$
as follows:
Let $d$ denote exterior differentiation, let $\ud 0\aa 0 = \aa 0$
and let $\ud 1\aa 0 = d\aa 0 \tmod \aa 0$.
 If $v_1\hd v_k\in T_f\cf^1$, extend
$v_1\hd v_k$ to holomorphic vector fields in some neighborhood of $f$
which we   denote $\tilde v_1\hd \tilde v_k$.
Let
$$
 \ud k\aa 0 (v_1\hd v_k ) :=  v_1(\intprod d (
\tilde v_2 \intprod
\hdots d(\tilde v_k\intprod d\aa 0 )
 \text{ mod }\pi_k\inv (\text{Image}\ud{k-1})
\tag 2.5
$$
where
$\pi_k: V\ra V/(\text{Image}\{\ud 0\ud 1\hd\ud{k-1}\} )$ is the projection,
and $\intprod$ denotes the contraction
$T\times  T^{*\ot l}\ra  T^{*\ot l-1}$. (2.5) is independent of
the extension of $v_1\hd v_k$ to holomorphic vector fields.
(The proof that (2.5) is independent of
 the choice of extension to holomorphic
vector fields is the same as the standard argument in the real
case, see e.g. [S].)

  $\ud k\aa 0$ descends to be a well defined element of
$$
  S^kT_x^*X\ot V/\tim
(\ud 0\hd \ud{k-1})\tag 2.6
$$
called
  the {\it   $k$-th fundamental form} of $X$ (twisted by $\Cal O(k-1)$)
which we will denote $\bii kX$, except that we will often denote
$\bii 2X$ by $II$ and $\bii 3X$ by $III$.

To fix notation, we will verify this assertion
in the case $k=2$. Given $f\in \cf^1$, let
$(\omega^B_C ) = \Omega = f\inv df$ denote the
entries of the Maurer-Cartan form.
Write $\oo B $ for $\ooo B 0$.
 The first two terms of   (2.5) expressed in frames are
 $$
\align
\ud 1\aa 0 &= \oo\alpha\ot\aa\alpha \text{ mod } \hat x\tag 2.7\\
\ud 2\aa 0 &= \oo\alpha\ooo\mu\alpha\ot\aa\mu\text{ mod}\hat T
 \tag 2.8
\endalign
$$
 To
see that $\ud 2\aa 0$ descends
(modulo twisting)
to a section of $S^2T^*\ot N$, note that for $g\in G_1$,
 $\Omega_{fg}=
g\inv\Omega_f g + g\inv dg$
(where the expression $dg$ is to  be understood as the
matrix of differentials
of the  functions $(g^A_B)$)
so tensoring (2.8)
with $\aa 0^*$;
the fiber motions  of $\oo\alpha, \ooo\mu\alpha,
\aa 0^*$ and $\aa\mu$ cancel.

Note that as a form on $\cf^1$,
$\oo\mu =0$ which implies $d\oo\mu=0$. (2.2) implies $d\oo\mu
=-\ooo\mu\alpha\ww\oo\alpha$. Since the forms $\{\oo\alpha\}$ are all
independent, the Cartan Lemma (see e.g. [BCG3]) implies that
$\ooo\mu\alpha =\qq\mu\alpha\beta\oo\beta$ for some functions
$\qq\mu\alpha\beta=\qq\mu\beta\alpha$ defined on $\cf^1$.
Thus (2.8) may be rewritten:
$$
\ud 2\aa 0  = \qq\mu\alpha\beta\oo\alpha\oo\beta \ot\aa\mu\text{ mod } \hat T.
$$

\smallpagebreak

For the purposes of this paper it will be useful to consider
$\bii 2X = II\in S^2T^*\ot N$ as a map
$$
II: N^*\ra S^2T^* \tag 2.9
$$
(which is dual to the standard way of considering $II$ as a map).

We will consider the  $k$-th fundamental form of $X$ at $x$    as a map
$$
\bii k{Xx} : \tker \bii{k-1}{Xx}\ra S^kT^*_x\tag 2.10
$$
where $\tker\bii{k-1}{Xx}\subset V^*$ and $\tker\bii 0{Xx}= \hat
x\upperp\subset V^* $.
  $\bpp{} (\tker \bii k{Xx} ) \subset \bpp{} V^*$  is the
 space of hyperplanes osculating to order $k$ at $x$.

\smallpagebreak

\noindent{\bf Remark 2.11}.
 We will often
omit reference to the base point $x$,
ignore twists, and use $T$ to denote both $\hat T/\hat x$ and $T_xX = \hat
x^*\ot \hat T/\hat x$.

\smallpagebreak

Another way of understanding fundamental forms is as follows:

The quotient  map
$$
  V^* \ra  V^*/ \hat x\upperp=\Cal O_{\bppp V }(1)_x \tag 2.12
$$
gives rise to a spectral sequence of a filtered complex by letting
$$
\align F^0K^0 &= V^*\ \ \  \ \ \  F^0K^1 = \Cal O_{X }(1)_x \tag 2.13 \\
F^1K^0 &= 0 \ \ \  \ \ \ \ \ F^p = F^pK^1 =   \frak m^p_x(1).
\endalign
$$
The maps are
$$
\align
& \ud 0 : V^*\ra F^0/F^1 = \Cal O_{X,x}(1)/\frak m_x(1)\simeq\bcc{}\tag 2.14
\\
& \ud 1: \text{ker}\ud 0  \ra F^1/F^2 = \frak m_x (1)/\frak m^2_x (1)
\simeq T^*(1) \\
& \ud 2: \text{ker}\ud 1 \ra F^2/F^3 = \frak m^2_x (1)/\frak m^3_x (1)
\simeq (S^2T^*)(1) \\
& \vdots \endalign
$$

\smallpagebreak

To study the space of hypersurfaces of
degree $d$ osculating
to order $k$ at a general point $x\in X$   in terms of   local differential
invariants,   we need to compute the  fundamental forms of
the $d$-th Veronese embedding of $X$,
$\vdx\subset\bpp{}S^dV$. We will do this in  \S 3, and their expression will
involve more subtle
invariants of $X\subset\bpp{}V$ which we now describe.

\smallpagebreak

Differentiating the equation $\ooo\mu\alpha-\qq\mu\alpha\beta\oo\beta=0$ and
using
the Cartan lemma, one obtains functions $\rr\mu\alpha\beta\gamma$ (symmetric in
their
lower indices) defined  on $\cf^1$  by the equation:
$$
\rr\mu\alpha\beta\gamma\oo\gamma =
-d\qq\mu\alpha\beta - \qq\mu\alpha\beta\ooo 00 -\qq\nu\alpha\beta\ooo\mu\nu
+\qq\mu\alpha\delta\ooo\delta\beta + \qq\mu\beta\delta\ooo\delta\alpha\tag 2.15
$$
(see [L1]). The form
$$
F_3 := \rr\mu\alpha\beta\gamma\oo\alpha\oo\beta\oo\gamma
\ot \aa\mu \text{mod }\hat T\tag 2.16
$$
 is a section of
$S^3(T^*\cf^1)\ot\pi^*(NX)$ (again, ignoring twists). $F_3$ was
defined for hypersurfaces in [GH] and  called
the {\it cubic form} there.  In [L1], $F_3$ was denoted $\partial II$.
We will also call $F_3$
{\it the first variation of II}. $F_3$
is actually a section of
$S^3(SB)\ot\pi^*NX$ where $SB\subset T^*\cf^1$ is the subbundle of
semi-basic forms, that is those that annhilate vertical tangent vectors.
 As defined, $F_3$ is a
{\it relative invariant} in the sense that it is defined as a section of a
bundle over
$\cf^1$ instead of a bundle over $X$. One can define $F_3$ as a section of a
bundle
over $X$ (see [L1]). However   it is not advantageous to do so because there is
a whole series of relative
invariants, of which $F_3$ is the first, and none of the others can be defined
as a section
of a bundle over $X$.  (Given a particular variety, one can make normalizations
that
enable one to define the relative invariants as a section of a bundle over $X$,
but the
normalizations
will not be canonical.) On the other hand, certain combinations of relative
invariants can be defined as sections of natural bundles over $X$. In this
paper we will deal with combinations that are the fundamental forms
of the
Veronese re-embeddings of $X$.

$F_3$ is the projective analogue of the covariant derivative of the second
fundamental
form  of a metric connection, which we will denote
$\nabla II$.     $\nabla II$ is defined as a section of a
bundle over the original variety. One may think of the projective structure as
specifying an equivalence class of connections and the necessity of defining
$F_3$ over a principal bundle as corresponding to the ambiguity in the choice
of
compatible connection.

Differentiating (2.15), one obtains functions $\rr\mu\alpha\beta{\gamma\delta}$
defined
on $\cf^1$  by
$$
\align
 \rr\mu\alpha\beta{\gamma\delta}\oo\delta &=
-d\rr\mu\alpha\beta\gamma -2\rr\mu\alpha\beta\gamma\ooo 0 0
-\rr\nu\alpha\beta\gamma\ooo\mu\nu + \frak
S_{\alpha\beta\gamma}\rr\mu\alpha\beta\delta\ooo\delta\gamma\tag 2.17 \\
& \ \ \
-\frak S_{\alpha\beta\gamma} \qq\mu\alpha\delta\qq\nu\beta\gamma\ooo\delta\nu +
\frak S_{\alpha\beta\gamma}\qq\mu\alpha\beta\ooo 0\gamma  \endalign
$$
 which leads to a form
$$
F_4 = \rr\mu\alpha\beta{\gamma\delta}\oo\alpha\oo\beta\oo\gamma\oo\delta
\ot \aa\mu \text{mod }\hat T\in\Gamma(S^4(T^*\cf^1)\ot\pi^*NX)\tag 2.18
$$

Say we have defined $F_{k-1}$ by differentiating the equations defining
the coefficients of $F_{k-2}$. Then we can define $F_k$
by differentiating the defining equations of the coefficients of
$F_{k-1}$.

We will call $F_k$ the {\it $(k-2)$-nd variation of $II$}.
   $F_k$ measures   how
$X$ is infinitesimally
leaving its embedded tangent space to order $(k-1)$ at $x$.
 We will   use the notation
$F_0$ for the quotient map $V^*\ra V^*/\hat x\upperp$,
  $F_1$ for the quotient map $\tker F_0\ra T^*$,
and  $F_2= \bii 2X$. On $\cf^1$,
$F_k$ is a section of
$S^kT^*\cf^1\ot\pi^*NX$, in fact a section of
$S^k(SB)\ot\pi^*NX$.

\proclaim{Proposition 2.19} The coefficients $\rr\mu{\alpha_1}\hdots
{\alpha_{l}}$ of $F_{l}$ are defined by the formula
$$
\align
\rr\mu{\alpha_1}\hdots
{\alpha_{l}}\oo{\alpha_{l}} &=\tag 2.20\\
&
-d\rr\mu{\alpha_1}\hdots{\alpha_{l-1}} -
l\rr\mu{\alpha_1}\hdots{\alpha_{l-1}} \ooo 0 0
-\rr\nu{\alpha_1}\hdots{\alpha_{l-1}}\ooo\mu\nu\\
&+
\frak S_{\alpha_1\hd\alpha_{l-1}}\rr\mu{\alpha_1}\hdots{\alpha_{l-2}\beta}
\ooo\beta{\alpha_{l-1}}
 +  l\frak S_{\alpha_1\hd\alpha_{l-1}}
\rr\mu{\alpha_1}\hdots{\alpha_{l-2}}\ooo 0{\alpha_{l-1}}\\
 &
-\Sigma_{p=1}^{l-3}\frak S_{\alpha_1\hd\alpha_{l-2}}
\rr\mu{\alpha_1}\hdots{\alpha_{p}\beta}
\rr\nu{\alpha_{p+1}}\hdots{\alpha_{l-2}}
\ooo \beta\nu
 \} \\
&
-\Sigma_{p=2}^{l-3}\frak S_{\alpha_1\hd\alpha_{l-1}} (p-2+l)
\rr\mu{\alpha_1}\hdots{\alpha_{p}}
\rr\nu{\alpha_{p+1}}\hdots{\alpha_{l-1}}
\ooo 0\nu .
\endalign
$$
\endproclaim
\demo{Proof}
Say we have computed the coefficients for $F_{k+1}$, i.e. that
$$
\align
0= &
d\rr\mu{\alpha_1}\hd{\alpha_k}-a_k\rr\mu{\alpha_1}\hd{\alpha_k}\ooo 00
-b_k\rr\nu{\alpha_1}\hd{\alpha_k}\ooo\mu\nu
+c_k\frak S_{\alpha_1\hd\alpha_{k}}
\rr\mu{\alpha_1}\hd{\alpha_{k-1}\beta}\ooo\beta{\alpha_k} \tag 2.21\\
&+e_k\frak S_{\alpha_1\hd\alpha_{k}}
\rr\mu{\alpha_1}\hd{\alpha_{k-1}}\ooo 0{\alpha_k}
+\Sigma_{p=1}^{k-2}f_{p,k}
\frak S_{\alpha_1\hd\alpha_{k}}
\rr\mu{\alpha_1}\hd{\alpha_{p}\beta}
\rr\nu{\alpha_{p+1}}\hd{\alpha_k}\ooo\beta\nu \\
&+
\Sigma_{p=2}^{k-2}g_{p,k}
\frak S_{\alpha_1\hd\alpha_{k}}
\rr\mu{\alpha_1}\hd{\alpha_{p} }
\rr\nu{\alpha_{p+1}}\hd{\alpha_k}\ooo 0\nu
-
\rr\mu{\alpha_{ 1}}\hd{\alpha_{k+1}}\oo{\alpha_{k+1}}
\endalign
$$
where $a_k\hd g_{p,k}$ are some constants independent of $X$.
Taking the exterior derivative of (2.21) the terms with semi-basic
coefficients (we can ignore the others since they   cancel out) are
$$
\align
&
a_k \rr\mu{\alpha_1}\hd{\alpha_{k+1}}\oo{\alpha_{k+1}}\ooo 00
+a_k\rr\mu{\alpha_1}\hd{\alpha_{k }}\ooo 0\beta\oo\beta\\
&+b_k\rr\nu{\alpha_1}\hd{\alpha_{k+1 }}\oo{\alpha_{k+1}}\ooo\mu\nu
+b_k\rr\nu{\alpha_1}\hd{\alpha_{k }}\ooo\mu\beta\ooo\beta\nu \\
&
+c_k\frak S_{\alpha_1\hd\alpha_k}
\{
\rr\mu{\alpha_1\hd\alpha_{k-1}}\beta{\alpha_{k+1}}\oo{\alpha_{k+1}}
 \ooo\beta{\alpha_k}
-\rr\mu{\alpha_1\hd\alpha_{k-1}}\beta{}
(\ooo\beta 0 \ooo 0{\alpha_k} + \ooo\beta\nu\ooo\nu{\alpha_k}) \}\\
&
-e_k\frak S_{\alpha_1\hd\alpha_k}
\{
\rr\mu{\alpha_1\hd\alpha_{k-1}}\beta{}\oo\beta\ooo 0{\alpha_k}
-
\rr\mu{\alpha_1}\hd{\alpha_{k-1}}\ooo 0\nu\ooo\nu{\alpha_k}\}\\
&
-\Sigma_{p=1}^{k-2}f_{p,k}
\frak S_{\alpha_1\hd\alpha_k}
\{
\rr\mu{\alpha_1}\hd{\alpha_{p}\beta\delta}
\rr\nu{\alpha_{p+1}}\hd{\alpha_{k}}\oo\delta\ooo\beta\nu
+
\rr\mu{\alpha_1}\hd{\alpha_{p}\beta}
\rr\nu{\alpha_{p+1}}\hd{\alpha_{k}\delta}\oo\delta\ooo\beta\nu \\
&\hphantom{-\Sigma_{p=1}^{k-2}f_{p,k}\frak S_{\alpha_1\hd\alpha_k}\{ }
  -\rr\mu{\alpha_1}\hd{\alpha_{p}\beta}
\rr\nu{\alpha_{p+1}}\hd{\alpha_{k}}\oo\beta\ooo 0\nu \} \\
&
-\Sigma_{p=2}^{k-1}g_{p,k}
\frak S_{\alpha_1\hd\alpha_k}
\{
\rr\mu{\alpha_1}\hd{\alpha_{p}\beta}
\rr\nu{\alpha_{p+1}}\hd{\alpha_{k}}
\oo\beta\ooo 0\nu
+
\rr\mu{\alpha_1}\hd{\alpha_{p} }\rr\nu{\alpha_{p+1}}\hd{\alpha_{k}\beta}
\oo\beta\ooo 0\nu \}\\
&
-d\rr\mu{\alpha_1}\hd{\alpha_{k+1} }\oo{\alpha_{k+1}}
+\rr\mu{\alpha_1}\hd{\alpha_{k+1} }(\oo{\alpha_{k+1}}\ooo 00
+\ooo{\alpha_{k+1}}\beta\oo\beta ).
\endalign
$$
 Collecting terms, the coeffeicient of $\oo{\alpha_{k+1}}$ is
$$
\align
&
-(a_k+1)\rr\mu{\alpha_1}\hd{\alpha_{k+1} }\ooo 00 \tag 2.22\\
& -
b_k\rr\nu{\alpha_1}\hd{\alpha_{k+1} }\ooo\mu\nu \\
&
+c_k\frak S_{\alpha_1\hd\alpha_k}
\rr\mu{\alpha_1}\hd{\alpha_{k+1 }\beta }\ooo\beta{\alpha_{k }}
+\rr\mu{\alpha_1}\hd{\alpha_{k }\beta }\ooo\beta{\alpha_{k+1}}\\
&
+a_k\rr\mu{\alpha_1}\hd{\alpha_{k } }\ooo 0{\alpha_{k+1}}
+c_k\frak S_{\alpha_1\hd\alpha_k}
\rr\mu{\alpha_1}\hd{\alpha_{k-1},\alpha_{k+1}}
\ooo 0{\alpha_k}
 +e_k
\frak S_{\alpha_1\hd\alpha_k}
\rr\mu{\alpha_1}\hd{\alpha_{k-1},\alpha_{k+1}}\ooo 0{\alpha_k} \\
&
+b_k\rr\nu{\alpha_1}\hd{\alpha_{k}}\qq\mu\beta{\alpha_{k+1}}\ooo\beta\nu
+c_k\frak S_{\alpha_1\hd\alpha_k}
\rr\mu{\alpha_1}\hd{\alpha_{k-1 },\beta}
\qq\nu{\alpha_k}{\alpha_{k+1}}\ooo\beta\nu \\
&
\Sigma f_{p,k} \frak S_{\alpha_1\hd\alpha_k}\{
\rr\mu{\alpha_1}\hd{\alpha_{p}\beta\alpha_{k+1} }
\rr\nu{\alpha_{p+1}}\hd{\alpha_k }\ooo\beta\nu +
\rr\mu{\alpha_1}\hd{\alpha_{p}\beta }
\rr\nu{\alpha_{p+1}}\hd{\alpha_k \alpha_{k+1} }\ooo\beta\nu \} \\
&
\frak S_{\alpha_1\hd\alpha_k}
\{ e_k\rr\mu{\alpha_1}\hd{\alpha_{k-1}}
\qq\nu{\alpha_k}{\alpha_{k+1}} \ooo 0\nu +
\Sigma f_{p,k}\rr\mu{\alpha_1}\hd{\alpha_p\alpha_{k +1}}
\rr\nu{\alpha_{p+1}}\hd{\alpha_{k}}\ooo 0\nu \\
&
+
\Sigma g_{p,k}\rr\mu{\alpha_1}\hd{\alpha_p\alpha_{k +1}}
\rr\nu{\alpha_{p+1}}\hd{\alpha_{k}}\ooo 0\nu \}.
\endalign
$$
Thus, comparing (2.22) with (2.15),(2.17) and (2.20),  we have
$a_{k+1}=a_k+1=k+1$, $b_{k+1}=b_k=1$, $c_{k+1}=c_k=1$,
$e_{k+1}= a_k =k$, $f_{p,k}=1$, $g_{p,k}=p-2+k$.
\qed\enddemo

To determine the actual geometric information in the invariants $F_k$ we must
know how they vary in the fiber.  $F_3, F_4$ vary in the fiber as follows:
If $(\tilde \aa 0,\tilde\aa\alpha,\tilde\aa\mu )$ is a new frame with
$$
\align
&\tilde\aa\mu = \aa\mu + g^0_{\mu}\aa 0 + g^{\alpha}_{\mu}\aa\alpha
\tag 2.23\\
& \tilde\aa\alpha = \aa\alpha + g^0_{\alpha}\aa 0\endalign
$$
then
$$
\align
&\tilde\rr\mu\alpha\beta\gamma = \rr\mu\alpha\beta\gamma +
\frak S_{\alpha\beta\gamma}g^0_{\alpha}\qq\mu\beta\gamma +
\frak S_{\alpha\beta\gamma}g^{\delta}_{\nu}\qq\nu\alpha\beta\qq\mu\gamma\delta
\tag 2.24\\
&\tilde\rr\mu\alpha\beta{\gamma\delta} =
\rr\mu\alpha\beta{\gamma\delta} +
\frak S_{\alpha\beta\gamma\delta}g^0_{\alpha}\rr\mu\beta\gamma\delta +
\frak S_{\alpha\beta\gamma\delta}g^{\epsilon}_{\nu}
(\rr\nu\alpha\beta\gamma\qq\mu\delta\epsilon+
\qq\nu\alpha\beta\rr\mu\gamma\delta\epsilon)
+ g^0_{\nu}\qq\mu\alpha\beta\qq\nu\gamma\delta .
\endalign
$$

We do not consider motions by $g^{\mu}_{\nu},g^{\alpha}_{\beta}, g^0_0$
because they just conjugate the coefficients by invertible matrices.
We will use the notation $\Delta\rr\mu\alpha\beta\gamma$ to denote
the change in $\rr\mu\alpha\beta\gamma$  by a  fiber motion of
the type in (2.23).
By (2.16),
$$
\Delta\rr\mu\alpha\beta\gamma =
\frak S_{\alpha\beta\gamma} (\gg 0\alpha\qq\mu\beta\gamma +
\gg\delta\nu\qq\nu\alpha\beta\qq\mu\gamma\delta ). \tag 2.25
$$

We will occasionally write
$F_k=\ff k\mu \aa\mu = \rr\mu{\alpha_1}\hdots{\alpha_k}\oo{\alpha_1}
\hdots\oo{\alpha_k}\aa\mu$ with the ambiguity understood.

Because they are only relative invariants, the invariants $F_k$ are more
difficult to deal with than fundamental
forms.

\smallpagebreak

\noindent{\bf Example 2.26}. Invariants of curves in $\bpp 2$.

(This computation is originally due to
 to Monge, and in this language   to Cartan [C].)
Take  an adapted frame $f=(\aa 0 ,\aa 1,\aa 2)$ defined up to the
action of $G_1$.
 As long as
our curve is not a line, at a general point we can rescale the
second fundamental form so that $\qq 211\equiv 1$. (e.g. by scaling
$\aa 0$).  We restrict to frames on which $\qq 211 \equiv 1$.
The defining equation for $F_3$ becomes
$$
  \rr 2111 \oo 1  = (-\ooo 0 0 -\ooo 22 + 2\ooo 11) \tag 2.27
$$
(ommiting the $\ot$ symbols here and in what follows) which we   normalize to
zero by sending
$$
\aa 1 \mapsto \aa 1 - \frac 12\rr 2111 \aa 0\tag 2.28
$$
(or by moving $\aa 2$ by $\aa 1$) and restricting to frames where
$\rr 2111 =0$.
Thus, for a curve in $\bpp 2$, $F_3$ contains no geometric information.
In our restricted frame bundle,
$$
  \rr 2{11}11\oo 1  = (\ooo 01 + \ooo 12) \tag 2.29
$$
which we normalize to zero by sending
$$
\aa 2 \mapsto \aa 2 - \rr 2{11}11\aa 0\tag 2.30
$$
so $F_4$ doesn't contain any geometric information either. Restricting to
frames where $\rr 2{11}11=0$,
$$
  \rr 2{111}11\oo 1  = \ooo 02. \tag 2.31
$$
$ \rr 2{111}11$ cannot be normalized to zero, but if it is not
zero, it can be made constant by scaling
$\aa 2$.  Making $\rr 211 {111}=1$   sets
$$
  \rr 2{1111}11\oo 1  = \ooo 01 \tag 2.32
$$
which  cannot be normalized, so $F_6$ represents the first  non-discrete
differential invariant of a nondegenerate curve in $\bpp 2$.
Note that $F_5=0$ implies $F_6=0$, and in fact that all higher $F_k=0$.
In \S 4 we will continue this example.

\smallpagebreak

\noindent{\bf Example 2.33.} Surfaces in $\bpp 5$ with zero third fundamental
form.

We may  adapt frames such that
$II_X= \oo1\oo 1\ot\aa 3 + \oo 1\oo 2\ot\aa 4 + \oo 2\oo 2\ot\aa 5\tmod\hat
T$.
The variability in $F_3$ is
$$
\matrix
 \Delta\rr 3111 = \gg 01 + \gg 13 & \qquad \Delta\rr 4111 =\gg 23 &\qquad
\Delta\rr 5111 = 0 \\
\Delta\rr 3112 = \gg 02 + \gg 14 & \qquad \Delta\rr 4112 =
\gg 01 + \gg 24 + \gg 13  &\qquad
\Delta\rr 5112 = \gg 21 \\
\Delta\rr 3122 = \gg 15   & \qquad \Delta\rr 4122 =
\gg 02 + \gg 14 + \gg 23  &\qquad
\Delta\rr 5122 = \gg 01 +\gg 24 \\
\Delta\rr 3222 = 0
&\qquad\Delta\rr 4222 = \gg 15
&\qquad
\Delta\rr 5222 = \gg 02 + \gg 25 \endmatrix \tag 2.34
$$
 Using the motions $\gg\alpha\mu$
in (2.23), we can normalize six of the twelve
terms $\rr\mu\alpha\beta\gamma$ to zero.  This uses
up all
our normalizations (the $\gg 0\alpha$ cannot contribute any additional
normalizations
in this case).
Informally,
we will say  $F_3$ contains $6$ functions worth of geometric information.

More generally, for $X^n\subset\bpp{n+\binom{n+1}2}$ with $III_X=0$, i.e.
$|II_X|=\bpp{}S^2T^*$, there are $\binom{n+1}2n$
effective normalizations $\gg\alpha\nu$ so
there are
 $\binom{n+1}2\binom{n+2}3
-\binom{n+1}2n$ functions worth of geometric information
in $F_3$.

\medpagebreak

For an $n$-dimensional variety of codimension $a$, if $F_3=0$, the
equations for  the coefficients of $F_4$ are
$$
\rr\mu\alpha\beta{\gamma\delta}\oo\delta =
\frak S_{\alpha\beta\gamma} (
\qq\mu\alpha\beta\ooo 0\gamma + \qq\mu\alpha\epsilon
\qq\nu\beta\gamma\ooo\epsilon\nu) \ \forall \mu,\alpha\leq\beta\leq\gamma\tag
2.35
$$
which is
$a\binom{n+2}3$ equations relating the $n$ independent forms $\{\oo\delta \}$
to the $n+an$ forms $\{\ooo 0\gamma , \ooo\epsilon\nu\}$.  (If one does not use
the normalizations
$\gg 0\alpha$, the forms $\{ \oo\alpha ,\ooo 0\beta \}$ are all independent and
the
$a\binom{n+2}3$ equations express the $an$ forms $\ooo\epsilon\nu$ in terms of
the $2n$ forms $\{ \oo\alpha ,\ooo 0\gamma \}$.)
Thinking of $II$ as given, the equations for the coefficients
of $F_3$
  are extremely overdetermined,   but sometimes
they have solutions.

\smallpagebreak

\noindent{\bf Example 2.36 } Griffiths and Harris showed that
if $n\geq 2$ and $a=1$, then
if one can restrict to   a sub-bundle $\cf^3\subset\cf^2$ on which the
functions
$\rr\mu\alpha\beta\gamma$ are all identically zero, then
one can restrict to a further subbundle $\cf^4\subseteq\cf^3$,
on which    the coefficients of
all the higher forms are also identically zero
([GH], B16).
To prove the $n=2$ case one must compute  the equations of $F_5$
to determine that the coefficients of $F_4 $ are zero.

\smallpagebreak

\noindent{\bf Example 2.37}. Let $X^n\subset\bpp{n+2}$ and
let  $x\in X$ be a general point.
Assume there is a nonsingular quadric in $\ii$. (This will be the case iff  $X$
has
a nondegenerate dual variety.   $X^*$ is nondegenerate if
$X$ is smooth and $n\geq 4$ by [E].)  We may normalize $II_X$ such that
$$
II_X =( (\oo 1)^2 + \hdots (\oo n)^2)\ot \aa{n+1} +
( \lambda_1(\oo 1)^2 + \hdots \lambda_n(\oo n)^2)\ot \aa{n+2}.\tag 2.38
$$
where $\lambda_i\in\bccc$.
This case is the closest to the hypersurface case   and
is still extremely overdetermined.  If $F_3=0$ we may let
$$
\rr\mu\alpha\beta{\gamma\delta} = \frak S_{\alpha\beta\gamma\delta}
\bb\mu\nu\tau\qq\nu\alpha\beta\qq\tau\gamma\delta     \tag 2.39
$$
where $\bb\mu\nu\tau =\bb\mu\tau\nu\in\bccc$ are constants.
Some, but not all of  these terms  can be normalized to zero.
(E.g. one can set $\bb\mu\mu\mu=0$ via $\gg 0\mu$.)
One can check that (2.39) is    a consistent solution to
the equations. i.e.

\smallpagebreak

{\it There exist varieties of
arbitrary dimension and codimension two such that it
is possible to reduce to frames where the coefficients of $F_3$ are
all identically zero,
   but that there is no frame such that  the coefficients of
$F_4$ are also  identically zero.}

\smallpagebreak

  We will
show that the equations above  are consistent in (4.34) with an easier,
indirect
proof.
The reader may wish to contrast this with the results of [JM].

\smallpagebreak

\noindent{\bf Remark 2.40}. In the metric situation, if
$\nabla II =0$  at a general point, then it is always the case that
$\nabla^kII=0$, $\forall k$ and the variety must be a locally
symmetric space. The reward of dealing with
our subtler invarients is that their  vanishing   will detect more general
properties than being a locally symmetric space.

 \smallpagebreak

Although example  (2.37) shows      that
being able to reduce to frames where the
coefficients of $F_k $ are all zero
does not imply that there is a further reduction
of frame such that the coefficients of  $F_{k+1} $ are also zero,
one does have  the following
result:

\proclaim{Proposition 2.41} Let $X\subset\ppp V$ be a variety, and
$x\in X$ any smooth point.
If there exists a framing where the coefficients
of $F_k\hd F_{2k} $ are all zero at $x$, then
the coefficients of
$F_l $ are also zero at $x$ $\forall l>2k$.
\endproclaim
\demo{Proof}
  Examining (2.20), we see that
if the coefficients
$\rr\mu{\alpha_1}\hd{\alpha_{k}}$
$\hd$
$\rr\mu{\alpha_1}\hd{\alpha_{2k }}$
are zero then the
$\rr\mu{\alpha_1}\hd{\alpha_{2k+1}}$
are as well.\qed
\enddemo

\bigpagebreak

\noindent\S 3.  {\bf Observations about the osculating hypersurfaces  of
  $X$}

As mentioned in the introduction,
$$\tker\bii k{\vdx x}= \{\text{ hypersurfaces of degree }
d\text{ osculating to order } k\text{ at }x\}.
$$
We have $\tker\bii k{\vdx x}\subseteq
\tker\bii {k+1}{\vdx x}$ and for sufficiently large $k$,
 the  maps
$\bii k\vdx$ are zero and the kernels
stabilize to $I_d\subset S^dV^*$.

In this section, we compute $\tker\bii k{\vdx x}$ in terms of the invariants
$F_l$
of $X$, make some observations about osculating hypersurfaces,
and describe a generalization of the
classical  Monge equation
that characterizes conic curves in $\pp 2$.

\smallpagebreak

The natural representation
$\rho_d : Gl(V) \ra Gl (S^dV)$ allows us to compute frames for
$\vdx$ using $\cf^1_{\bppp V}$ instead of
$\cf^1_{\bppp (S^dV)}$.
More precisely, one may compute the fundamental forms of
$\vdx$ using any adapted frame bundle over $\vdx$ that respects
the filtration $\hat x\subset\hat T_x\vdx\subset S^dV^*$. In
particular, we may restrict to the frame bundle
$\rho_d (\cf^1_X)\subset\cf^1_{\vdx}$.
 We   compute the
$k$-th fundamental form of
$\vdx$ by computing $\ud k (\aa 0^d)$.

In other words, we compute
the fundamental forms of $\vdx$ by applying the operator $\ud{}$ to
  sections of
  $\Cal O (d) \ra X$.
In the spectral sequence perspective,   one  uses
the  quotient  maps
$$
  S^dV^*\ra S^dV^*/(\hat x\upperp )^d=\Cal O_{X }(d)_x \tag 3.1
$$
and the analogous filtrations.
The equality
$
\Cal O_X(d) = v^*_d (\Cal O_{v_d(X)}(1)) $ connects the perspectives.

\smallpagebreak

\noindent{\bf Example 3.2}. Fundamental forms of $\vtwox$.

Assume $II_X$ is injective, write $\aaa 02$ for $\aa 0\aa 0$.  The
fundamental
forms of $\vtwox$ expressed in frames are:
$$
\align  \ud 1\aaa 0 2 &= 2\oo\alpha\aa\alpha\aa 0
\ \text{mod} \{ \aaa 0 2 \} \tag 3.3 \\
 \ud2\aaa 0 2 &= 2\oo\alpha\oo\beta(\qq\mu\alpha\beta\aa\mu\aa 0
+\aa\alpha\aa\beta)
\ \text{mod} \{ \aa\alpha\aa 0 ,\aaa 0 2 \}  \\
\ud 3\aaa 0 2 &=  2\oo\alpha\oo\beta\oo\gamma (
\rr\mu\alpha\beta\gamma \aa\mu\aa 0 + 3\frak S_{\alpha\beta\gamma}
\qq\mu\alpha\beta\aa\mu\aa\gamma) \\
&\ \ \ \text{mod} \{\qq\mu\alpha\beta\aa\mu\aa 0 +\aa\alpha\aa\beta,
 \aa\alpha\aa 0 ,\aaa 0 2 \}
\\
\ud 4\aaa 0 2 &=  2\oo\alpha\oo\beta\oo\gamma\oo\delta
(\rr\mu{\alpha\beta}\gamma\delta\aa\mu \aa 0 +
4\frak S_{\alpha\beta\gamma\delta}\rr\mu\alpha\beta\gamma\aa\mu\aa\delta
+ 3\frak S_{\alpha\beta\gamma\delta}\qq\mu\alpha\beta\qq\nu\gamma\delta
\aa\mu\aa\nu ) \\
 &\  \text{ mod }\{\text{Image}(\ud 3, \ud 2 , \ud 1, \ud 0)\}\\
&\vdots
  \endalign
$$
or in invariant notation,
$$
\align
\bii 1 \vtwox &= 2F_1F_0|_{\hat x^2\upperp} \tag 3.4\\
\bii 2 \vtwox &= 2(F_2F_0 + F_1F_1)|_{ (\hat x\hat T)\upperp}  \\
\bii 3 \vtwox &= 2(F_3F_0 + 3F_2F_1) |_{\tker\bii 2\vtwox} \\
\bii 4\vtwox&= 2(F_4F_0 + 4F_3F_1 + 3F_2F_2)|_{\tker\bii 3\vtwox}\\
 \bii 5\vtwox&= 2(F_5F_0 + 5F_4F_1 + 10F_3F_2)|_{\tker\bii 4\vtwox}\\
&\vdots
\endalign
$$

\noindent{\bf Remark 3.5}.
If we had not assumed that
$II_{Xx}$ was injective,   then the coefficients of
$III_X$ would have appeared beginning in $\ud 3$ and
its
infinitesimal variations
   would have appeared in higher order terms.

\smallpagebreak

\noindent{\bf 3.6. The classical Monge equation, Example 2.26 continued}.
Examining (3.4), we see that $F_5=0$ implies $\bii 5{v_2(C)} =0$.
If $x$ is a general point of $C$, this implies
$I_2=\tker\bii 4\vtwox$.
On the other hand, $\tdim (\tker\bii j{v_2(C)})=5-j$ for $j\leq4$ and therefore
there is a $Q\in S^2V^*$ osculating to all orders, i.e.   $C$ is a conic curve.
   $F_5=0$ is exactly the classical Monge equation for
a curve in a plane to be a conic.

 In fact, if we work in local affine coordinates $[1,x,y]$ where the
curve is $y=y(x)$ we can recover the original
Monge equation.  First take a lift
 $$
f=(\aa 0 ,\aa 1,\aa 2 )=
\pmatrix 1&0&0\\ x&1&0 \\ y&0&1\endpmatrix
$$
and then solve for $g\in G_1$ such that
$$
(fg)\inv d(fg)=\pmatrix \ooo 00 & \ooo 01 & \ooo 02\\
\ooo 10 & \ooo 11 & \ooo 12\\\ooo 20 & \ooo 21 & \ooo 22\endpmatrix
\tag 3.7
$$
satisfies
$\ooo 20  =\ooo 21 -\ooo 10=\ooo 12+\ooo 01=
\ooo 00 + \ooo 22 -2\ooo 11 =0$. Let $\tilde f=fg$. In our lift,
$\ooo 10=dx$ and $\ooo 02= (*)((y'')^{-\frac 23})'''dx$,
where $(*)$ is a nonzero term of lower order. Setting $\ooo 02=0$
yields the Monge equation
$$
 ((y'')^{-\frac 23})'''=0. \tag 3.8
$$
 The computation above is a little involved, and is
substantially simplified if one instead lifts such that
$\oo 1=(y'')^{-\frac 32}dx$ (this makes
$\oo 1$ have unit \lq\lq affine arc-length\rq\rq ). The computation
in affine space is carried out in detail in [C] and [S], Vol. II.

\smallpagebreak

\noindent{\bf 3.9 Example 2.36 continued}.
The actual Griffiths-Harris  result mentioned in
(2.36) is that for $n>1$, hypersurfaces  with $F_3=0$ are
quadrics.

 \smallpagebreak

\proclaim{Proposition 3.10}
The fundamental forms of $v_d(X)$ are
$$
\bii k \vdx =
\Sigma_{l_1+\hdots +l_d=k}c_{l_1\hdots l_d}F_{l_1}\hdots
F_{l_d} \tmod (\Sigma_{l<k} \bii l\vdx (S^lT) )|_{\tker\bii {k-1}\vdx} \tag
3.11
$$
where the $c_{l_1\hdots l_d}$ are nonzero constants.
\endproclaim

\demo{Proof}
The two main facts are
$$
\underline d^r (A\circ B) = \Sigma_{a+b=r}(\underline d^a A)
\circ (\underline d^b B)\ \  \text{ (Leibinitz rule)} \tag 3.12
$$
and
$$
\underline d^r \aa 0 = F_r. \tag 3.13
$$
Thus computing $\underline d^k (\aa 0 )^d$ is just like computing the $k$-th
derivative of a function of one variable raised to the $d$-th power.
I.e.
 let $f(x)$ be a function of one variable. Then
$$
(\frac d{dx})^kf^d(x)=(f^d)^{(k)}=
\Sigma_{l_1+\hd l_d=k}c_{l_1\hdots l_d}f^{(l_1)}\hdots f^{(l_d)}
\tag 3.14
$$
In particular, all the
$c_{l_1\hdots l_d}$ with $l_1+\hdots + l_d=k$ are nonzero.
One way of computing the coefficients is to use induction via the
formula
$$
(f^d)^{(r)}
 = \Sigma_{k=0}^r\binom r k f^{(r-k)}(f^{(d-1)})^{(k)}.\qed
\tag 3.15
$$
\enddemo

\smallpagebreak

Note that as long as $p\leq d$, $|\bii p{\vdx}|_x=\ppp (S^pT^*)$
for all $x\in X$
because of the $F_1^pF_0^{d-p}$ term in (3.11).
This property  is unchanged if the higher fundamental forms of $X\subset\bpp{}
V$
are nonzero.  This observation
is a generalization of the statement
 that
at smooth points $x\in X$ there is always
exactly  an $(a-1)$ dimensional family of hyperplanes osculating
to
first order. It implies

 \proclaim{Proposition  3.16}
 Let $X^n\subseteq\bpp{} V=\Bbb C\pp{n+a}$ be a
variety and let $x\in \xsm$.
 For all $  p\leq d$,

$$
\align &
\tdim
\left\{
\matrix  \text{(not necessarily irreducible)
hypersurfaces} \\\text{of degree }d\text{ osculating   to order }  p  \text{
at }x\endmatrix \right\} \\
&
= \binom{(n+a+1)+(d-1)}  d -  \{1+
n+ \binom{ n+1} 2  + \hdots +
\binom{ n+p-1}p
 \}.\endalign
$$
\endproclaim
\demo{Proof}
$$
\align \tdim (\tker &\bii p{v_d(X) x})\\
& =
(\tdim(S^dV^*)- \Sigma_{k=0}^{p-1}\tdim(\tim\bii k\vdx ))
 - \tdim(\tim\bii p\vdx )
.\qed
\endalign
$$
\enddemo

The pattern  in (3.16) does not continue. For example, it is not true that
the space of hypersurfaces of degree $d$ osculating to order $d+1$ is of
dimension max
$ \{ 0, \binom{(n+a+1)+(d-1)}  d -  \{
n+ \binom{ n+1} 2  + \hdots +
\binom{ n+(d+1)-1}{d+1}
 \}\}$. The actual dimensions of the spaces of osculating hypersurfaces depend
on
properties of the differential invariants of $X$ that we will discuss shortly.
However, independent
of the structure of the invariants, there are
{\it lower}
 bounds on the dimensions
for
$d+1\leq k\leq 2d-1$.

To see this, first note that the filtration
$\hat x\subset \hat T\subset V$ induces a
  $(2d+1)$-step  filtration of
$S^dV$.  For example, the seven step filtration of
$S^3V$ induced is
$$
\{\hat x^3\}\subset
\{\hat\xx 2\hat T \}\subset\{\hat\xx 2V,\hat x\hat T^2 \}
\subset \{\hat x\hat T V,\hat T^3 \}\subset
\{\hat x V^2,\hat T^2V\}
\subset \{\hat T V^2 \}\subset \{ V^3\}
$$
where  $\hat \xx 3=\hat x\circ\hat x\circ\hat x,\hat x\hat T=\hat x\circ\hat T$
etc...

Considering the fundamental forms dually as maps
$$
{\bii {*k}\vdx} : S^kT^*\ra N_{\vdx}/\tim\bii{k-1}\vdx ,
$$ the image
of  $\bii{*k}\vdx$ lies in a quotient of the $(k+1)$-st term
of the filtration of $S^dV$. Since the filtration
has $(2d+1)$ steps, it is not possible
for $\tker\bii k\vdx$ to be zero for $k\leq 2d$.
This leads to {\it lower} bounds on $\tdim (\tker\bii k\vdx)$
for $d<k\leq 2d-1$. For example:

 \proclaim{Proposition  3.17} Let $X^n\subseteq\bpp{} V=\Bbb C\pp{n+a}$ be any
variety,
 and $x\in X $, any smooth point.
$$
\tdim
\left\{
\matrix  \text{(not necessarily irreducible )
hypersurfaces of}\\\text{degree }d\text{  osculating   to order } 2d-1  \text{
at }x\endmatrix \right\}
  \geq
\binom {a+d-1}d-1.
$$
 \endproclaim
\demo{Proof}
The first term in any $\bii k\vdx$ for which $S^dN^*$ is not in the kernel
  is
$(F_2)^d$, which appears in $\bii {2d}\vdx$. Thus
$S^dN^{*}\nsubseteq\tker\bii{2d-1}\vdx$. Finally note, that
 dim$S^dN^{*}=\binom {a+d-1}d$.\qed
\enddemo

 \smallpagebreak

Another way to phrase the upper bounds on the dimensions of
the space of osculating hypersurfaces is in terms of lower
bounds on the dimensions of the osculating spaces
For example,
$$
\align
& \tdim \tilde T_x\up 2\vtwox = n+\binom{n+1}2 \tag 3.18\\
& \tdim \tilde T_x\up 3\vtwox \leq n+\binom{n+1}2 + a(1+n) \\
& \tdim \tilde T_x\up 4\vtwox \leq n+\binom{n+1}2 + a(1+n) +\binom{a+1}2\\
& \hphantom{\tdim \tilde T_x\up 4\vtwox}\ =  \tdim \ppp S^2V\endalign
$$
 and
$$
\align
& \tdim \tilde T_x\up 3\vthreex = n+\binom{n+1}2+\binom{n+2}3  \\
& \tdim \tilde T_x\up 4 \vthreex\leq n+\binom{n+1}2+\binom{n+2}3 + a(1+n) \\
& \tdim \tilde T_x\up 5 \vthreex\leq n+\binom{n+1}2
+\binom{n+2}3+ a(1+n) +\binom{a+1}3\\
& \tdim \tilde T_x\up 6\vthreex \leq n+\binom{n+1}2
+\binom{n+2}3+ a(1+n) +\binom{a+1}3 +\binom{a+2}3\\
& \hphantom{\tdim \tilde T_x\up 4\vthreex}\ =  \tdim \ppp S^3V.
\endalign
$$

We  now rephrase (1.6) in a manner that  will allow us to approximate it
using projective differential invariants and interpret it geometrically.

There is a natural codimension $a$ subspace of $N^*_{\vdx}$,
namely $(N\hat \xx{d-1})\upperp$, which in particular
contains all
the hypersurfaces in $I_d$ that are singular at $x$.

\proclaim{Proposition 3.19} Let $x\in X\subset \bppp V$ be a smooth point.
$X$ is a complete intersection if and only if
for all $k$,
$$
I_k\cap  (N\hat \xx{k-1})\upperp
\equiv 0 \ \tmod \ (I_{k-1}\circ V^*)\cap (N\hat \xx{k-1})\upperp .
$$
\endproclaim
\demo{Proof}
By  (1.6), it is sufficient to show that
$$
  \tker [d]_{  k  }  =I_k\cap  (N\hat \xx{k-1})\upperp
  \tmod   (I_{k-1}\circ V^*)\cap  (N\hat \xx{k-1})\upperp . \tag 3.20
$$

Let $(\xx 0,\xx\alpha ,\xx\mu )$ be a basis of $V^*$ (which we also think of as
linear coordinates on $V$)
such that $[1,0\hd 0] =x$, and $\hat T = \{ \xx\mu=0 \}$. Then  $N_x^*$ is
spanned
by $\{ d\xx\mu\}$.
Expressed in terms of our basis, elements of $(N\hat x^{k-1})\upperp$
have no $\xx\mu(\xx 0 )^{k-1}$ terms.
If $P\in I_k$ is such that $dP_x\neq 0$ then
$P= y(\xx 0 )^{k-1} + \hdots$
where $y=c_{\mu}\xx\mu$
for some constants $c_{\mu}$, and therefore $P\notin (N\hat x^{k-1})\upperp$.
If $[dP_x]=0$, but $dP_x\neq 0$ then $P=\xx 0 P'+P''$ where $P'\in I_{k-1}$
and $dP_x''=0$.\qed
\enddemo

One can define conditions that imply
the condition $(CI)_d$ by requiring that for some fixed $k$,
the analog of (3.19) holds for the hypersurfaces of degree
$d$ osculating to order $k$. More precisely, that
 $$
  \tker \bii{k}\vdx\cap { (N\hat \xx{d-1})\upperp}\equiv 0
\tmod  (I_{d-1}\circ V^*)\cap { (N\hat \xx{d-1})\upperp}.\tag 3.21.k
$$

\smallpagebreak

\noindent{\bf Definition 3.22}. Let $X\subseteq\ppp V$ be a variety and
$x\in X_{sm}$. We will say $X$
{\it is known to have no excess equations
of  degree $d$ after taking  $k$
derivatives at $x$},
or that {\it $X$ satisfies $(CI)^k_d$ at $x$}, if (3.21.k) holds
at $x$.

 We will work with the  condition $(CI)^{2d}_d$. By (3.10), $k=2d$ is
the smallest possible $k$ for which  (3.21.k) could hold.
We have no reason to believe this
is the \lq\lq correct\rq\rq\ condition but it provides a useful starting point.
(By \lq\lq correct\rq\rq\  we mean  that it will be possible to prove
some global theorem that  forces the condition to hold under suitable
codimension and smoothness hypotheses
about $X$.)  We will require these restrictions hold
at {\it general points} instead of all smooth points. This will enable us to
compute
using only a fixed, finite number of derivatives.

 \proclaim{Proposition 3.23 (generalized Monge systems) }
Fix positive integers $d_1<d_2<\hdots <d_r$ and
$a_1,a_2,\hdots , a_r$, where
$a_1+\hdots + a_r=a$.

 If $X^n\subset\bppp V=\Bbb C\pp{n+a}$
is a complete intersection and
$x\in X$ is a general point such that for all $d_j$
and for all $k$ such that $d_{j-1}<k\leq d_j$ (set $d_0=0$),
$$
\tdim\{  \tker \bii{k }{v_{k}(X)} \} \tmod (N\hat\xx{d_j-1})\upperp
=
a_1+\hdots + a_j \tag 3.24
$$
and
$$(CI)^{2k}_{k} \text{ holds at }x \tag 3.25
$$
then $X$ is a complete intersection of $a_1$ hypersurfaces of degree
$d_1$, $a_2$ hypersurfaces of degree $d_2$, ... ,
$a_r$ hypersurfaces of degree $d_r$.

Moreover,
$
\{ dP_x |  P\in \tker \bii{2d_j }{v_{d_j}(X)} \}
= N_{d_j}^*
$
where $N_{d_j}^*$ is the $j$-th term in the natural filtration
of $N^*_xX$ (1.4).
\endproclaim
\demo{Proof}
The assumption of
$(CI)^{2k}_k$ at $x$ implies that the only possible elements
of $\tker\bii{2k-1}{v_k(X),x}$ correspond to smooth hypersurfaces.
(3.24)
implies
$\bii{2d_j+1}{\vdx , x} =0$
so $I_d=\tker\bii {2d_j}{\vdx ,x}$. Furthermore,
(3.24) implies
  $\tdim (I_{d_j}/(I_{d_j-1}\circ
V^*))= a_j$.\qed\enddemo

Proposition (3.23) needs some explanation. (3.24) may be understood
as systems of partial differential equations. In fact,
if one
writes out the differential invariants in terms of local coordinates
then (3.24) is a system of PDE of order $2d_r+1$ which specializes to
the classical Monge equation (3.6)   when $n=a=1$ and $d=2$.

The   conditions $(CI)^{2d_j}_{d_j}$ expressed in terms of local coordinates
are also conditions on the
derivatives, but
at least in small codimension,
 they are open conditions that one might hope to
prove must  be satisfied in certain geometric situations.
The proposition implies:
\smallpagebreak

{\it Among
complete intersections that satisfy suitable genericity
conditions, if $d$ is the largest degree of a hypersurface essentially
containing $X$, the the ideal of $X$ can be recovered by taking
$2d+1$ derivatives at a general point.}
\smallpagebreak

In practice, it is actually useful to
 work with a stronger condition than
$(CI)^{2d}_d$. To understand the stronger condition,
interpret
$(CI)^{2d}_d$ as saying that $\tker\bii{2d}\vdx$ is as small as
possible
modulo the differentials of the polynomials
that are smooth at $x$. The stronger condition will be
for $d<k\leq 2d$,
that
  $\tker\bii{k}\vdx$ be as small as possible
modulo the differentials of the polynomials
that are smooth at $x$.
In this paper we will only discuss this stronger condition in the
case $d=2$ and will explain it in the next section.

 \bigpagebreak

\noindent\S 4. {\bf Osculating quadrics}

\smallpagebreak

We now discuss the   quadrics containing $X\subset\bppp V$ in detail.
We  assume throughout this section that $III_X=0$ at general points.

\smallpagebreak

Consider the fundamental forms of $\vtwox$
at a point $x$ restricted to the
conormal directions corresponding to the osculating quadrics
singular at $x$.
 (We use the notation
$T=\hat T/\hat x, N=V/\hat T$):
$$
\align
&
\bii 1\vtwox : T^*\hat\xx*\ra T^* \text{
\hphantom{xxxxxxxx}(Identity on
the first factor)}   \tag 4.1
\\
&
\bii 2\vtwox |_{(N\hat x )\upperp} :
S^2T^*  \ra S^2T^* \ \ (  \text{Identity map)}\\
&
\bii 3\vtwox |_{(N\hat x )\upperp} :
N^*T^*
\ra  S^3T^* \ \  (  II\circ\text{Identity map} )\\
&
\bii 4\vtwox |_{(N\hat x )\upperp} :
(\tker\bii 3\vtwox |_{(N\hat x )\upperp}\cap N^*T^*) +   N^*N^*
\ra  S^4T^* \\
&\ \ \hphantom{xxxxxxxxxxxxx} ( F_3\circ \text{ Identity map } + II\circ II )
\endalign
$$

$\tker \bii 3\vtwox |_{(N\hat x )\upperp}$ will be as small as possible if
$$
F_2F_1 :    (NT)^* \ra S^3T^*\tag 4.2
$$
is injective and assuming
(4.2) is injective,
 $\tker\bii 4\vtwox |_{(N\hat x )\upperp}$  will be empty if
$$
F_2F_2 :  N^{2*} \ra S^4T^*\tag 4.3
$$
is also injective.

We will say that
$X$ satisfies
{\it strong genericity
in degree two } at $x$ if  (4.2) and (4.3) are
injective.

One can similarly define strong genericity in degree $d$
at $x$, but it is more complicated
because one needs to take
$I_{d-1}$ into account.

It will be useful to study (4.2), (4.3) in bases.
Let $(\xx 0,\xx \alpha ,\xx\mu )$ be an adapted  basis of $V^*$,
dual  to  $(\aa 0,\aa\alpha,\aa\mu)$.
We also think of  $(\xx 0,\xx \alpha ,\xx\mu )$ as linear coordinates on $V$.

The induced basis of $S^2V^*$ is $(\xx 0\xx 0,\xx\alpha \xx 0,
\xx\mu\xx 0,\xx\alpha\xx\beta ,\xx\mu\xx\beta ,\xx\mu\xx\nu )$.
(3.4) implies
  $$
\align
\tker \bii 0\vtwox &= \{\xx\alpha \xx 0,
\xx\mu\xx 0,\xx\alpha\xx\beta ,\xx\mu\xx\beta ,\xx\mu\xx\nu \}
\tag 4.4 \\
\tker\bii 1\vtwox &=
\{\xx\mu\xx 0,\xx\alpha\xx\beta ,\xx\mu\xx\beta ,\xx\mu\xx\nu \}\\
\tker\bii 2\vtwox &=
\{ \xx\mu\xx 0 -
\qq\mu\alpha\beta\xx\alpha\xx\beta,\xx\mu\xx\beta,\xx\mu\xx\nu\}
\endalign
$$

 Beyond this, the   the spaces of osculating quadrics
depend on the structure of the invariants of $X$.

$\bii 3\vtwox$ maps its domain as follows:
$$
\align
\xx\mu\xx 0 - \qq\mu\alpha\beta\xx\alpha\xx\beta&\mapsto
2\ff 3\mu
\tag 4.5\\
\xx\mu\xx\gamma&\mapsto
6\ff 2\mu\oo\gamma\\
\xx\mu\xx\nu&\mapsto 0
\endalign
$$
where we are writing $\ff 3\mu = \rr\mu\alpha\beta\gamma\oo\alpha\oo\beta
\oo\gamma$ etc ...  (Here, since we have chosen a particular frame, the
coefficients $\rr\mu\alpha\beta\gamma$ etc... are well defined).

The failure of  (4.2)  to be injective  means that there exists
a nontrivial  equation $l_{\mu}q^{\mu}=0$ with $l_{\mu}\in T^*$,
$q^{\mu}\in\widehat\ii$.  We will call such an equation a
{\it linear syzygy} in $\ii$.

Similarly,
(4.3)  will fail to be injective if there exists a nontrivial equation
$
k_{\mu\nu}q^{\mu}q^{\nu}=0
$
with $k_{\mu\nu}=k_{\nu\mu}\in \Bbb C$.
We will refer to such an equation as a
{\it quadratic relation among the quadrics in} $\ii$.

In fact,
requiring (4.3) to be injective  is redundant:

\proclaim{Lemma 4.6}
Let $T$ be a vector space. If a system of quadrics
$A\subseteq S^2T^*$   satisfies
a quadratic relation among the quadrics in $A$, i.e. an equation of the form
$
c_{\mu\nu}q^{\mu}q^{\nu}=0
$
with $q^{\mu}$ a collection of independent elements of $A$ and
$c_{\mu\nu}=c_{\nu\mu}$ constants (not all zero), then $A$ has linear syzygies.
\endproclaim
\demo{Proof}
Let $v\in T$. Then
$$
v\intprod (c_{\mu\nu}q^{\mu}q^{\nu})= 0
$$
i.e.
$$
c_{\mu\nu}(v\intprod q^{\mu})q^{\nu}=0
$$
which is a nontrivial linear syzygy
with $l_{\nu} = c_{\mu\nu}(v\intprod q^{\mu})$
as long as $v\notin\text{Singloc}(
c_{\mu\nu}q^{\mu})$ for some $\nu$. Since not all
the quadrics $c_{\mu\nu}q^{\mu}$ are zero, we can always find such a $v$.
\qed
\enddemo

Summarizing:

\proclaim{Proposition 4.7} Let $X\subset\bpp{} V$ be a  variety.
 Let $x\in X$
be a general point and assume $III_{Xx}=0$ and
$II_{Xx}$ and has no linear syzygies.  Then
$X$ satisfies $(CI)_2$ at $x$, i.e. $X$ has no excess equations in
degree $2$
at $x$.\endproclaim

\noindent{\bf Remark 4.8.}
One can compare (4.7) with the following two varieties:
 The twisted cubic curve
$C\subset\bpp 3$  is cut out by quadrics and is not a complete
 intersection,
but  $III_C\neq 0$. The Segre $X=\bpp 1\times\bpp 2\subset\bpp 5$, or
any generic hyperplane section of it, is also cut out by quadrics
and not a complete intersection, but  there is a linear syzygy among
the quadrics in $II_X$ at any point due to the
point in the  base locus
(see 5.4 for more details).

\smallpagebreak

We now derive
 a more refined version of (3.23) for intersections of quadrics.

In order that $N_x^*$ be spanned by differentials of quadratic
polynomials, it is necessary that
$$
\{ dP_x | P\in \tker \bii k\vtwox \} = N_x^*\tag 4.9.k
$$
for all $k$.
(We supress reference to the base
point $x$ in what follows.)
 For $k\leq 2$, (4.9.k) automatically holds; for
$k=3$ (4.9.3) will hold if and only if
$$
\ff 3\mu = 3\al\mu\nu\gamma\oo\gamma\ff 2\nu
\tag 4.10
$$
for some constants $\al\mu\nu\gamma\in \Bbb C$.
Notice that if
$\rr\mu\alpha\beta\gamma =  \frak S_{\alpha\beta\gamma}
\al\mu\nu\gamma \qq\nu\alpha\beta  $
in some frame, it holds in any choice of frame (with different
constants $\al\mu\nu\gamma$), so the expression (4.10) has intrinsic meaning.
If (4.8) holds, then
 $$
\tker\bii 3\vtwox =
\{\xx\mu\xx 0 - \qq\mu\alpha\beta\xx\alpha\xx\beta-\al\mu\nu\beta\xx\nu
\xx\beta,
\xx\mu\xx\nu \}.\tag 4.11
$$
Assuming (4.10), and that there are no linear syzygies in $II_X$,
  $\bii 4\vtwox$  maps its domain as follows:
$$
\align
\xx\mu\xx 0 - \qq\mu\alpha\beta\xx\alpha\xx\beta
-\al\mu\nu\alpha\xx\nu\xx\alpha
&\mapsto
 2\ff 4\mu - 8\al\mu\nu\alpha\oo\alpha\ff 3\nu\tag 4.12\\
\xx\mu\xx\nu
&\mapsto
6\ff 2\mu\ff 2\nu
\endalign
$$
(4.9.4) is the condition
$$
\ff 4\mu = 4\al\mu\nu\alpha\oo\alpha\ff 3\nu + 3 \bb\mu\nu\tau
\ff 2\nu\ff 2\tau\tag 4.13
$$
for some constants $\bb\mu\nu\tau =\bb\mu\tau\nu\in\bccc$
(this expression also has intrinsic meaning, in that if it holds in
one choice of frame, it will hold in all choices).  If (4.10),
(4.13) hold and there
 are no linear syzygies in $II_X$, then
$$
\tker\bii 4\vtwox = \{ \xx\mu\xx 0-\qq\mu\alpha\beta\xx\alpha\xx\beta -
\al\mu\nu\gamma\xx\nu\xx\gamma -
\bb\mu\nu\tau\xx\nu\xx\tau
\}.
 \tag 4.14
$$
  $\bii 5\vtwox$  maps its domain as follows as follows:
$$
\xx\mu\xx 0-\qq\mu\alpha\beta\xx\alpha\xx\beta -
\al\mu\nu\gamma\xx\nu\xx\gamma -
\bb\mu\nu\tau\xx\nu\xx\tau
\mapsto
2\ff 5\mu - 10\al\mu\nu\gamma\oo\gamma\ff 4\nu - 20\bb\mu\nu\tau\ff 3\nu\ff
2\tau
\tag 4.15
$$
(4.9.5) is the condition
$$
\ff 5\mu = 5\al\mu\nu\gamma\oo\gamma\ff 4\nu + 10\bb\mu\nu\tau\ff 3\nu\ff
2\tau\tag 4.16
$$
(which also has intrinsic meaning).

 If (4.10),(4.13),(4.16) and strong genericity all hold, then
  $\bii 5\vtwox =0$. Since we are at a general
point,  this implies all higher
fundamental forms are zero and   $I_2= \tker\bii 4\vtwox$.

In summary:
 if
$$
\align &
 \ff 3\mu = 3\al\mu\nu\gamma\oo\gamma\ff 2\nu \tag 4.17\\
&\ff 4\mu = 4\al\mu\nu\alpha\oo\alpha\ff 3\nu + 3 \bb\mu\nu\tau
\ff 2\nu\ff 2\tau\\
&\ff 5\mu = 5\al\mu\nu\gamma\oo\gamma\ff 4\nu + 10\bb\mu\nu\tau\ff 3\nu\ff
2\tau
\endalign
$$
  where $ \al\mu\nu\alpha , \bb\mu\nu\tau=\bb\mu\tau\nu\in \Bbb C$,
and strong
genericity   in degree $2$ holds at $x$,
then
 $N_x^*$ is  spanned by the differentials of a set of generators
of  $I_2$. I.e. the only hypersurfaces
that essentially contain $X$ are of degree two. In this
case,   we will call (4.17) the
{\it generalized Monge system for quadrics}.

 In summary:

\proclaim{Theorem 4.18}
Let $X\subset\bppp V$ be a variety and $x\in X$ a general point.
Assume $III_{Xx}=0$   and that there are
no linear syzygies in $ \ii_x$.
Then
$$
\align &
\tdim \{ \text{quadrics osculating to order three at } x \}\leq
a
+\binom{a+1}2-1  \tag 4.19\\
&
\tdim \{ \text {quadrics osculating to order four at } x \}\leq
a-1.
\endalign
$$

If  the generalized Monge system (4.17) holds, then
$$
 I_2=\tker\bii 4{\vtwox x}  =
\{ \xx\mu\xx 0-\qq\mu\alpha\beta\xx\alpha\xx\beta -
\al\mu\nu\gamma\xx\nu\xx\gamma -
\bb\mu\nu\tau\xx\nu\xx\tau,  \ n+1\leq\mu\leq n+a
\}
$$
where $\qq\mu\alpha\beta$ are the coefficients of the second fundamental form
at $x$ and $\al\mu\nu\alpha ,\bb\mu\nu\tau$ are   coefficients
expressing $F_3,F_4$ in terms of $F_2$.

    Equality occurs in the first (respectively second)  line of  (4.19)
 if and only if    the first
(resp. second) line of (4.17) holds at $x$.    If the
generalized Monge system  does not hold, then $I_X$ is not generated
 by quadrics.
 \endproclaim

Theorem (4.18) implies:
\smallpagebreak

 \proclaim{Corollary 4.20} Let $X\subset\bppp V$ be a variety and $x\in X$ a
general point.
Assume $III_{Xx}=0$. If   there are no quadric hypersurfaces singular at $x$
that osculate to order four at $x$, and $Q$ is a quadric hypersurface
osculating to order five at $x$, then $X\subseteq Q$.\endproclaim

\smallpagebreak

\noindent{\bf Remark 4.21.} The assumption of small codimension (which is
implicit in
the hypotheses that $\bii 3{Xx}=0$   and that there are
no linear syzygies in $ \ii_x$) is essential to being able to determine the
quadratic equations of $X$ by taking only five derivatives.  One already needs
more derivatives for a curve in $\bpp 3$.

   \proclaim{Corollary 4.22 }  $X$ as in (4.18)  is determined by
$II_X,F_3, F_4$ at one general point. In fact, the
higher variations of $II$ are given by the formula:
$$
c_{k0}\ff k\mu  =\frac 12 ( c_{ {(k-1)}1}\al\mu\nu\alpha\oo\alpha\ff{k-1}\nu
+\bb\mu\nu\tau (\Sigma _{l+m=k} c_{ lm}\ff l\nu\ff m\tau ))  \ \forall  k>2.
\tag 4.23
$$
 \endproclaim

\smallpagebreak

\noindent{\bf Example 4.24.}
Let $X^6\subset\bpp {12}$,  let  $x\in X$  be a smooth point  and let
$$
II_{Xx} = \oo 1\oo 4\aa 7 + \oo 2\oo 5\aa 8 + \oo 3\oo 6\aa 9 +
\oo 1\oo 5\aa {10} + \oo 2\oo 6\aa {11}+\oo 3\oo 4\aa{12}.\tag 4.25
$$
Say in addition that  $F_3=F_4=0$.  Then $I_X$ is generated by$$
\align
 &
\{\xx 0\xx 7 - \xx 1\xx 4 ,\xx 0\xx 8 - \xx 2\xx 5,
\xx 0 \xx 9 - \xx 3\xx 6 , \xx 0 \xx {10} - \xx 1\xx 5 ,
\xx 0 \xx{11} -\xx 2\xx 6,\tag 4.26\\
& \
\xx 0 \xx{12} -\xx 3\xx 4,
\xx 7\xx 8\xx 9 - \xx{10}\xx{11}\xx{12} \}.\endalign
$$

\smallpagebreak

The presence of    relations or syzygies among the quadrics in $\ii$
can produce  equations of higher degrees
that are not elements of
$I_2\circ S^kV^*$
as in (4.26).  Not all syzygies and relations
actually produce such equations, there are  tautological ones, which we will
call the {\it Koszul syzygies}.
For example, let $\ii = \{ Q^{\mu}\}$, where $Q^{\mu} =
\qq\mu\alpha\beta\oo\alpha\oo\beta$.
For each $\mu <\nu$ there are Koszul syzygies
$$
\qq\nu\alpha\beta\oo\alpha \oo\beta Q^{\mu}
-\qq\mu\alpha\beta\oo\alpha \oo\beta Q^{\nu}=0.\tag 4.27
$$

\smallpagebreak

\proclaim{Theorem  4.28}
Let $X\subset\ppp V$ be a variety.
A sufficient condition for   $I_X$ to be generated by   quadrics
is that at a smooth point $x\in X$:

1.
$\{ dP_x | P\in I_2 \} = N^*_x$

2. Any  sygyzies  or polynimials satisfied by $\ii_x$ of the form
$$
l_{\mu_1\hdots\mu_k ,\alpha_1\hdots\alpha_p}q^{\mu_1}\hdots
q^{\mu_k}\oo{\alpha_1\hdots\alpha_p}=0
$$
where $l_{\mu_1\hdots\mu_k ,\alpha_1\hdots\alpha_p}\in\bccc$,  other than
the  {\it Koszul} syzygies,   are generated by the linear syzygies
and the quadratic relations among the quadrics in $\ii_x$.

In particular, if
$X\subset\ppp V$ is a variety, $x\in X$ is  a general point,
there are no linear syzygies in $\ii_x$,
no polynomials satisfied by the quadrics in $\ii_x$, and the generalized Monge
system
 (4.17) holds,
then   $I_X$ is generated by quadrics, in fact
$X$ is a complete intersection of
quadrics.
 \endproclaim
\demo{Proof}
Let $(\xx 0,\xx\alpha ,\xx\mu )$ be an adapted basis of $V^*$
 dual to $(\aa 0 ,\aa\alpha ,\aa\mu )$.

Consider the case of cubics containing $X$.
Condition 1. implies that $I_2$ contains quadrics of the
form
$$
\xx\mu\xx 0 +\hdots \tag 4.29
$$
for each $n+1\leq \mu\leq n+a$ so
$I_2\circ V^*$ contains cubics of the form
$$
x^B(\xx\mu\xx 0 +\hdots) \tag 4.30
$$
$\forall 0\leq B\leq n+a$.  Given $P\in I_3$,
we may modify  $P$ by elements of $I_2\circ V^*$ so that we
may assume
$$
P(\aa{B\mu 0})=0 \ \forall 0\leq B\leq n+a. \tag 4.31
$$
That $P$ vanishes to all orders on $X$ says $P(\bii k\vthreex )=0 \ \forall k$.
(Where we are considering the pairing $  S^3V^*\times
\bii k\vthreex \ra S^kT^*$.)
In particular,

$$
\align
&P(\bii 0\vthreex )=0\Rightarrow P(\aa{000})=0 \tag 4.32
\\
&P(\bii 1\vthreex )=0\Rightarrow
P(\aa{\alpha 00})=0\ \forall\alpha
\\
&P(\bii 2\vthreex )=0\Rightarrow
P(\aa{\alpha\beta 0})=0\ \forall\alpha,\beta
\\
&P(\bii 3\vthreex )=0\Rightarrow
P(\aa{\alpha\beta\gamma})=0\ \forall\alpha,\beta,\gamma
\\
&P(\bii 4\vthreex )=0\Rightarrow
P(\frak S_{\alpha\beta\gamma\delta}\qq\mu\alpha\beta\aa{\mu\gamma\delta})
= \frak S_{\alpha\beta\gamma\delta}\qq\mu\alpha\beta
P(\aa{\mu\gamma\delta})
=0\
\forall\alpha,\beta,\gamma,\delta
\\
&P(\bii 5\vthreex )=0\Rightarrow
P(
\frak S_{\alpha\beta\gamma\delta\epsilon}
\qq\mu\alpha\beta\qq\nu\gamma\delta\aa{\mu\nu\epsilon})=0\
\forall\alpha,\beta,\gamma,\delta,\epsilon
\\
&P(\bii 6\vthreex )=0\Rightarrow
P( \frak
S_{{\alpha_1}\hdots{\alpha_6}}\qq\mu{\alpha_1}{\alpha_2}
\qq\nu{\alpha_3}{\alpha_4}\qq\tau{\alpha_5}{\alpha_6}\aa{\mu\nu\tau})=0\
\forall {\alpha_1}\hd{\alpha_6}.\endalign
$$
where in each line we have used the line above and (4.31)
to reduce to only having one nonzero term to worry about.
Using (4.32) and the absense of polynomials and linear syzygies,
 we see  that the modified $P$ is zero.
i.e. that
$I_3\equiv 0 \tmod I_2\circ V^*$.

The proof for equations of degree $d$ is the same. One may assume
that $P$ contains no terms of the form $\xx{B_1}\hdots\xx{B_{d-2}}\xx\mu\xx 0$
and   use the first $2d+1$ fundamental forms of $\vdx$.
\qed\enddemo

\medpagebreak

 We now explain what happens in the special case of (4.17) when
$F_3=F_4=0$:

\proclaim{Theorem 4.33}Let $X\subset \bppp V$ be a variety and
$x\in X$ a general point. Assume
$III_{Xx}=0$. If $F_3=F_4=0$ at $x$, then
the  minimal
  number of generators of $I_X$ is  $a+$
the minimal number of
generators of
     syzygies  and polynomials
in $\ii_x$ modulo the Koszul syzygies.  Furthermore,
either $X$ is a smooth homogeneous
space or $X_{sing}\subset H$, a hyperplane, and $X\backslash
(X\cap H)$ is a homogeneous
variety in the affine space $\bppp V\backslash H$. In either
case $X$ is birationally equivalent to $\bpp n$ and
$I_X$ is generated by quadrics if
and only if all the syzygies  and polynomials of $ \ii$ are generated by the
linear
  syzygies and quadratic relations among the quadrics in $\ii$ and the
 Koszul  syzygies and relations.
\endproclaim

\demo{Proof}
$X$ is given by the following construction (see [LVdV], [L1]):

Fix any smooth point $x\in X$.
Let $Y= \text{Baseloc}\ii_x\subset\bppp T_x$.
 Linearly embed $T\subset\Bbb C^{n+1}$.
Let $y^0\hd y^n$ be
linear coordinates on $\Bbb C^{n+1}$ such that $T= \{ y^0=0\}$.
Consider the rational map:
$$
Bl_Y\bpp n - - -\ra \bppp \{\widehat{\ii} , y^0\circ T^*\}^*\subseteq
\bppp \{ \Cal O_{\bpp n}\ot \ci_Y\} \subset\bppp (S^2\Bbb C^{n+1}) \tag 4.34
$$
To see that the image of this map is isomorphic to $X$,
 note that it has the correct codimension and that  all differential
invariants at $x$ are the same (only $II_X$ is nonzero for both $X$ and the
image), thus the varieties are
isomorphic.\qed
\enddemo

 \smallpagebreak

\noindent{\bf Remark 4.35.} If $X$ is as in (4.33),
let $ \bpp{a-1} = \{ \xx 0 = 0 ,\xx\alpha =0\}$.
$X$ will be singular in the
space
 $X\cap \bpp{a-1}$ unless the  number
 of generators of
   syzygies and polynomials modulo Koszul is
at least $a-1$. In fact, if there are no polynomials among the
quadrics, then $\bpp{a-1}\subset X$ (see (4.39)).

\smallpagebreak

\noindent{\bf 4.36. Example 2.37 continued}.
(2.39) may be rephrased as $F_4 = \bb\mu\nu\tau \ff 2\nu\ff 2\tau$.
After normalizing  $\bb {n+1}{(n+1)}{(n+1)}= \bb {n+2}{(n+2)}{(n+2)}=0$, all
our effective
normalizations that keep $F_3=0$ are used up and
the equations of $X$ are
$$
\align
&
\{ \xx {n+1}\xx 0 - \Sigma_{\alpha}\xx\alpha\xx\alpha -
b^{n+1}_{(n+1)(n+2)}\xx{n+1}\xx{n+2}-b^{n+1}_{(n+2)(n+2)}\xx{n+2}\xx{n+2}, \\
&\ \
\xx {n+2}\xx 0 - \Sigma_{\alpha}\lambda_{\alpha}\xx\alpha\xx\alpha -
b^{n+2}_{(n+1)(n+1)}\xx{n+1}\xx{n+1}-b^{n+2}_{(n+1)(n+2)}\xx{n+1}\xx{n+2}
\}
\endalign
$$
 Note that
since it is   impossible to normalize
both the $\al\mu\nu\alpha$'s and the $\bb\mu\nu\tau$'s to zero simultaneously,
$X$ cannot be rational.

\smallpagebreak

\noindent{\bf Example 4.37}.
We describe the equations of    varieties  as in (4.33)
whose ideals are generated by quadrics.
Write the map (4.34) as
$$
x\mapsto [y^0y^0, y^0y^1 \hdots y^{0}y^n, Q^1(x)\hd Q^a(x)]
=[x^0\hd x^n\hd x^{n+a}]
$$
Say there are $r$ linear syzygies
$k^l_{\mu\alpha}\xx\alpha Q^{\mu}$, $1\leq l\leq r$ and
$s$  quadratic relations among the quadrics in  $\ii$,
$k^m_{\mu\nu}Q^{\mu}Q^{\nu}$, $1\leq m\leq s$
and assume they and the  Koszul  syzygies generate the space
of syzygies and relations.   The $a+r+s$ equations of
$X$ are
$$
\{
x^0x^{\mu}-\qq\mu\alpha\beta x^{\alpha} x^{\beta},
\, k^l_{\mu\alpha}x^{\mu}x^{\alpha},\, k^m_{\mu\nu}x^{\mu}x^{\nu}\}.
\tag 4.38
$$

\noindent {\bf 4.39  An apparent dichotomy for varieties cut out by quadrics}.

It is extremely difficult
for a variety $X\subset\bppp V$
with $III_{Xx}=0$ at general points $x\in X$, and
satisfying (4.17) to have extra equations
without $X$ being rational.  For example, say there  is a linear syzygy of the
form
$$
k_{\mu\alpha}\oo\alpha\ff 2\mu=0 \tag 4.40
$$
where $k_{\mu\alpha}\in\bccc$.  In order for this linear syzygy to produce an
extra  equation, the constants $\al\mu\nu\alpha ,\bb\mu\nu\tau$ must satisfy
the equations
$$
\align
&
k_{\mu\alpha}\oo\alpha\al\mu\nu\beta\oo\beta\ff 2\nu=0 \tag 4.41\\
&
k_{\mu\alpha}\oo\alpha(\al\mu\nu\beta\oo\beta\al\nu\tau\gamma\oo\gamma\ff 2\tau
+\bb\mu\nu\tau\ff 2\nu\ff 2\tau) =0\\
&
k_{\mu\alpha}\oo\alpha(\al\mu\nu\beta\oo\beta
(\al\nu\tau\gamma\oo\gamma \al\tau\sigma\delta\oo\delta\ff 2\sigma +
\bb\nu\tau\sigma\ff 2\tau\ff 2\sigma)
+\bb\mu\nu\tau\al\nu\sigma\delta\oo\delta\ff 2\sigma\ff  2 \tau) =0\endalign
$$
which are severely overdetermined, on the order of
$n^5$ equations for the $a^2n+\binom{a+1}2a$ coefficients
$\al\mu\nu\alpha ,\bb\mu\nu\tau$.

\proclaim{Question 4.42} Let $X$ be a variety and let $x\in X$ be a general
point.
If $I_X$ is generated by quadrics    and
$\bii 3{Xx}=0$,
  is it necessarily the case that  $X=Z_1\cap Z_2$
where $Z_1$ is  a complete intersection and $Z_2$ is
 rational, and both $I_{Z_1}$, $I_{Z_2}$ are generated by quadrics?
\endproclaim

\noindent{\bf Remark 4.43.} Note that if $III_X\neq 0$ the equations (4.41)
are replaced
with a much less overdetermined system, so the  validity  of the question
is  heavily dependent
on $III_X=0$.   If $a\leq \frac{n-(b+1)}2$, then we are guaranteed
 $III_X=0$ by ([L1], (4.15)).

\smallpagebreak

We will show the answer to   question (4.42) is {\it yes}
 in the codimension range
 $a<\frac{n-(b+1) +2}3$ in \S 6.  If the answer to (4.42) is no, then (4.41)
gives a guide to all potential counter-examples to Hartshorne's conjecture
on complete intersections
whose ideals  are generated by quadrics.  The first
intersesting case is for $11$-folds in $\bpp {15}$.

The system one must solve to get a non-complete intersection is as follows:

Let $1\leq i,j,k\leq r$, $r+1\leq s,t\leq n$.  One must find constants
$\qqq ijk ,\qqq ijs ,\al{n+i}{n+j,}k$,
$\al{n+i}{n+j,}s , \bb{n+i}{n+j,}{n+k}$ as follows:
Let
$\rr{n+i}\alpha\beta\gamma=\frak S_{\alpha\beta\gamma}
\al{n+i}\mu\alpha\qq\mu\beta\gamma$ etc... be given as in (4.17).
One needs
$$
\align
&\frak S_{ijk}\qqq ijk =0, \ \qqq ijs + \qqq jis =0\tag 4.44\\
&\frak S_{ijkl}\rr{n+i}jkl = 0, \
\frak S_{ijk}\rr{n+i}jks = 0,  \  \rr{n+i}jst +\rr{n+j}ist =0\endalign
$$
and the analogous equations for $F_4$ and $F_5$ to hold.
If one is looking for smooth non-complete intersections, one also
needs to check smoothness, which amounts to genericity conditions on
$\al\mu\nu\alpha,\bb\mu\nu\tau$.

 \proclaim{Proposition 4.45}  The answer to (4.42) is yes for surfaces in $\bpp
4 $
and  $\bpp 5$ and $3$-folds in
$\bpp 5$.
\endproclaim

\demo{Proof}
Case of surfaces in $\bpp 4$.
$\ii$ must consist of two quadrics if $III_X=0$. In order for  $X$ to fail
to be a complete intersection, there must be a linear syzygy in $\ii$. This
implies
we can choose frames such that  $II_X = \oo 1\oo 1\ot \aa 3 + \oo 1\oo 2 \ot
\aa 4$
in a neighborhood of our general point $x$. The variablility in the
coefficients of $F_3$ is as follows
$$
\matrix
\Delta\rr 3111 = 3\gg 13 + \gg 01 &\qquad\qquad\qquad\qquad\qquad
 \Delta\rr 4111 = 3\gg 23
 \\
  \Delta\rr 3112 = 2\gg 14 + \gg 02
&\qquad\qquad\qquad\qquad\qquad
 \Delta\rr 4112 = \gg 13 + \gg 24 + \gg 01\\
 \Delta\rr 3122 = 0 & \qquad\qquad\qquad\qquad\qquad\qquad
 \Delta\rr 4122 = \gg 14 + \gg 02\\
 \Delta\rr 3222 = 0 &\qquad\qquad\qquad\qquad\qquad\qquad
 \Delta\rr 4222 = 0
\endmatrix\tag 4.46
$$
Using $\gg\alpha\mu$ set
$\rr 3111 ,\rr 3112 ,\rr 4111 , \rr 4112 =0$.  If  $I_X$ is  generated by
quadrics, (4.17) implies
$$
\align &
\rr 3111 = 3\al 331 ,\ \rr 3112 = \al 332 + 2\al 341 ,\ \rr 3122 = 2\al 342
,\rr 3222 = 0
\tag  4.47
\\
&
\rr 4111 = 3\al  431 ,\ \rr 3112 = \al 432 + 2\al 341 ,\ \rr 4122 = 2\al 442
,\rr 3222 = 0
\endalign
$$
which implies $ \al 331 = \al 332 + 2\al 341 =
\al 342 = \al 431 =\al 432 + 2\al 441 =0$ in these frames.

To have the syzygy persist, the relation
$ \oo 1\ff 34 - \oo 2\ff 33 =0$ must also hold,  which
forces all the $\al\mu\nu\alpha =0$ and thus $F_3=0$.

Similarly,  if $I_X$ is generated by quadrics,
$\ff 4\mu = \bb\mu\nu\tau\ff 2\nu\ff 2\tau$ which,
after normalizing by $\gg 03 ,\gg 04$ and
requiring $\oo 1\ff 44 - \oo 2\ff 43=0$ implies $F_4=0$.

The case  of a $3$-fold in $\bpp 5$ is similar,
only  there are two possibilities for $\ii$,
namely $\{ \oo 1\oo 1 ,\oo 1\oo 3 \}$ and $\{ \oo 1\oo 2 , \oo 1\oo 3 \}$.
The case
of a surface in $\bpp 5$ was proved in [L2] (assuming even less
than we  assumed here).\qed
\enddemo

\noindent\S 5. {\bf Some Homogeneous Examples}

In this section we write out the equations of some homogeneous varieties in a
manner
that illuminates the computations of \S 4.

\smallpagebreak

\noindent {\bf Example 5.1,  The Veronese}. $v_2(\bpp{}W)\ra\bpp{}S^2W$.

Let $W$ have basis $(B_0,B_{\alpha})$, $1\leq\alpha,\beta\leq n$.
Let $\aa 0 = B_0\circ B_0,\aa\alpha=B_0\circ B_{\alpha}$.
Let $\mu = (\alpha , \beta ),\ \alpha\leq\beta$ index the normal directions,
i.e. $\aa\mu = B_{\alpha}\circ B_{\beta}$.
If $\gamma\leq \delta$ we may adapt frames such that
$$II_X=\oo\alpha\oo\beta\aa{(\alpha\beta )}
\tmod \{\aa 0,\aa\alpha \}\tag 5.2
$$
I.e.,
 $\qq{(\alpha , \beta )}\gamma\delta
= \delta^{\alpha}_{\gamma}\delta^{\beta}_{\delta}$. The equations are
 $$
\{ \xx 0 \xx{(\alpha , \beta )}-\xx\alpha\xx\beta ,
\xx{(\alpha , \beta )}\xx\gamma - \xx{(\alpha , \gamma )}\xx\beta
 - \xx{(\beta , \gamma )}\xx\alpha,
\xx{(\alpha , \beta )}\xx{(\gamma , \delta )} -
\xx{(\alpha , \gamma )}\xx{(\beta , \delta )} \}\subset S^2V^* \tag 5.3
$$
(where in the second term we require $\alpha <\beta <\gamma$)
of which there are
$a+\binom n 3+  [\binom n4 + 3\binom n3 ]=\tdim\Lambda^2(S^2W)$.
  The first set of equations  comes from
$\ii$, the second linear syzygies in   $ \ii$, and the third  quadratic
relations among the quadrics in  $ \ii$.

\smallpagebreak

\noindent{\bf Example 5.4, The Segre}. $\bpp{}W_1\times\bpp{}W_2 \ra
\bpp{}(W_1\ot W_2 )$.

Let $(B_0, B_{\alpha})$, $(C_0, C_j )$ be respective adapted framings
of $W_1,W_2$, $1\leq i,j\leq m$, $1\leq\alpha,\beta\leq n$.  Let $\aa 0 =
B_0\ot C_0$, $\aa\alpha = B_{\alpha}\ot C_0$,
$\aa j = C_j\ot B_0$ and $\aa{\alpha j} = B_{\alpha}\ot C_j$.
So $\{\aa 0\}\subset\{\aa 0, \aa\alpha,\aa j\}\subset\{\aa 0,
\aa\alpha,\aa j,\aa{\alpha j}\}$ is a first order adapted framing.
Letting $\{\oo\alpha\}$,
$\{\phi^i\}$ be the pullbacks of the semi basic forms on $\bppp W_1$
and  $\bppp W_2$  respectively, then
$$
II_X = \oo\alpha\phi^j \ot  \aa{\alpha j} \text{ mod } \{\aa 0, \aa\alpha ,
\aa j \}.\tag 5.5
$$
The equations  of the Segre  are
$$
\{\xx 0\xx{\alpha j} - \xx\alpha\xx j, \
\xx j\xx{\alpha k} - \xx k \xx{\alpha j} ,\
\xx\beta\xx{\alpha k}- \xx\alpha\xx{\beta k} , \
\xx{\alpha j}\xx{\beta k} - \xx{\beta j}\xx{\alpha k} \}\tag 5.6
$$
where the first set of equations comes from $II_X$, the second and
third from  linear syzygies in $ \ii$, and the last set
from  quadratic relations among the quadrics in  $ \ii$, a total of
$a + \left[ \binom n2 m + \binom m 2 n \right]
 + \binom n2\binom m 2=(\tdim\Lambda^2V)(\tdim\Lambda^2W )$ equations.

\smallpagebreak

\noindent{\bf Example 5.7, The Grassmannian} $G(2, W)\subset
\bpp{}(\Lambda^2W)$.

Identify the tangent space to $G(2,W)$ with the $2\times (m-2)$ matrices
(dim$W=m$) and index everything accordingly.
Write the normal indices as $\mu = (ij) , i<j, 3\leq i,j\leq m$
and the tangent indicies as $ \alpha = (1j)$ or $(2j)$.
Then
$$
II_X =  (\oo{1j}\oo{2k} - \oo{2j}\oo{1k})\ot  \aa{jk}
\text{ mod } \{\aa 0, \aa\alpha\}.\tag 5.8
$$
The equations are
$$
\align \{ &\xx 0 \xx{ij} - (\xx{1i}\xx{2j} - \xx{1j}\xx{2i}) ,\
\xx{1k}\xx{ij} + \xx{1i}\xx{jk} - \xx{1j}\xx{ik},\tag 5.9\\
 & \xx{2k}\xx{ij} + \xx{2i}\xx{jk} - \xx{2j}\xx{ik}, \
\xx{ij}\xx{kl} + \xx{il}\xx{jk} -\xx{ik}\xx{jl} \}\endalign
$$
where the first set of equations comes from $II_X$, the second and
third from linear syzygies in $ \ii$, (which comes from picking two
columns
and a tangent vector) and the last set
from  quadratic relations among the quadrics in  $ \ii$ (which comes
from picking
pairs out of four columns), a total of
$a + 2\binom m 3 + \binom m 4=\tdim\Lambda ^4 W$ equations.

\smallpagebreak

\noindent{\bf Example 5.10, The Severi Varieties}. (In particular,
$E_6/P\subseteq\bpp{26}$).

Refer to [L2] for notations.  Let $\Bbb A$ be a division algebra over $\Bbb C$.
Let  $\Bbb A^2$
have division algebra valued coordinates $(u,v)$.  For $X$ Severi,
 $T_xX\simeq \Bbb A^2$ and
$$
\ii = \bppp\{ u\overline v , u\overline u , v\overline v \}\tag 5.11
$$
where the
first expression gives dim$_{\Bbb C}\Bbb A$ quadrics and the other
two expressions give one each.  The  linear sygyzies are
$$
\overline u (u\overline v) - \overline v (u\overline u ), \ \
v(u\overline v) - u (v\overline v) \tag 5.12
$$
where each expression gives dim$_{\Bbb C}\Bbb A$ linear relations.
There is a unique quadratic relation among the quadrics in $\ii$;
$$
(u\overline v )(\overline u v) - (u\overline u ) (v\overline v) \tag 5.13
$$
where $\overline u v$ is the same set of quadrics as
$u\overline v$.
In fact we can see how these fit into the set of all equations  by
writing an element of $y\in \Cal H = V$ as
$$
y = \pmatrix r_1 & \bar u_1 & \bar u_2 \\
u_1 & r_2 & \bar u_3 \\ u_2 & u_3 & r_3 \endpmatrix \ \ r_i\in \Bbb C , \
u_i\in \Bbb A \tag 5.14
$$
and taking
$$
\aa 0 = \pmatrix 1& 0&0\\ 0&0&0\\0&0&0\endpmatrix .
$$
The equations for the Severi variety are just the $2\times 2$ minors.
Here $r_1$ plays the role of $x^0$ in the previous examples.  The equations
are:
$$
\{r_1r_2-u_1\overline{u_1} , r_1r_3 - u_2\overline{u_2},
r_1u_3-\overline{u_1}u_2 , r_2u_2 - u_1u_3, r_3u_1-u_2\overline{u_3}, r_2r_3 -
u_3\overline{u_3} \}\tag 5.15
$$
where the first three terms ($2+\tdim_{\Bbb C}\Bbb A$ equations) come from
$\ii$, and using them to rewrite the
rest without $r_1,r_2, u_3$, we see the fourth and fifth terms
($2\tdim_{\Bbb C}\Bbb A$ equations) are the linear syzygies in
$ \ii$ and the last term (1 equation) is the
quadratic  relation among the quadrics in  $ \ii$.
The reader may find it amusing to explicitly correlate this description of
$v_2(\bpp 2)$, Segre$(\bpp 2\times\bpp 2)$, and $G(2,6)$ with the ones given
above.

\smallpagebreak

\noindent{\bf Example 5.16. The Spinor variety} $\Bbb S^{10}\subset\pp{15}$.

Let $V=\bcc 5$.
Write $\bcc{16}=\Lambda^{even}V=\Lambda^0V\oplus \Lambda^2V\oplus\Lambda^4V$
with basis
$(\aa 0,\aa{ij},\aa i )$, $1\leq i<j\leq 5$. Fixing $x\in\Bbb S$,
we may identify $x\simeq [\aa 0]\simeq\bppp
(\Lambda^0 V)$, $\hat T_x\Bbb S\simeq
\{ \aa 0,\aa{ij} \}\simeq (\Lambda ^0V\oplus \Lambda ^2 V)$ and the system of
quadrics obtained from
the second fundamental form is the complete system of quadrics
with base locus the Grassmanian $G(2,5)\subset\ppp (\Lambda^2V)$.
We may take
$$
\align II_X\equiv &
(\oo{12}\oo{34}-\oo{13}\oo{24}-\oo{14}\oo{23})\ot\aa 5\tag  5.17\\
&+(\oo{12}\oo{35}-\oo{13}\oo{25}-\oo{15}\oo{23})\ot\aa 4\\
&+(\oo{12}\oo{45}-\oo{14}\oo{25}-\oo{15}\oo{24})\ot\aa
3\\
&+(\oo{13}\oo{45}-\oo{14}\oo{35}-\oo{15}\oo{34})\ot\aa
2\\
&+(\oo{23}\oo{45}-\oo{24}\oo{35}-\oo{25}\oo{34})\ot\aa 1
\tmod \{\aa 0 ,\aa{ij} \}.\endalign
$$
Here there are only
linear syzygies. The equations are
$$
\align
 \{ & \xx 5\xx 0 - (\xx{12}\xx{34}-\xx{13}\xx{24}-\xx{14}\xx{23}),\
\xx 4\xx 0 -(\xx{12}\xx{35}-\xx{13}\xx{25}-\xx{15}\xx{23})\tag 5.18\\
&
 \xx 3\xx 0-(\xx{12}\xx{45}-\xx{14}\xx{25}-\xx{15}\xx{24}),\
\xx 2\xx 0 -(\xx{13}\xx{45}-\xx{14}\xx{35}-\xx{15}\xx{34}),\\
&  \xx 1\xx 0-(\xx{23}\xx{45}-\xx{24}\xx{35}-\xx{25}\xx{34}) ,\\
&\xx{15}\xx 5 + \xx{14}\xx 4+\xx{13}\xx 3 + \xx{12}\xx 2,\
\xx{25}\xx 5 + \xx{24}\xx 4+\xx{23}\xx 3 - \xx{12}\xx 1,\\
&\xx{35}\xx 5 + \xx{34}\xx 4-\xx{23}\xx 3 - \xx{13}\xx 1,\
\xx{45}\xx 5 - \xx{34}\xx 3-\xx{24}\xx 2 - \xx{14}\xx 1,\\
&\xx{45}\xx 4 + \xx{35}\xx 3+\xx {25}\xx 2 +\xx{15}\xx 1\}.
\endalign
$$

The reader may find it amusing to explicitly correlate
these equations with the system of quadrics (5.11) with
$\Bbb A=\Bbb O$.

\medpagebreak

 \noindent\S 6.  {\bf Properties of Systems of Quadrics}

\smallpagebreak

In this section we study  what the existence of
a linear syzygy
implies about a system of quadrics
and  the implications of a companion condition. We then combine these
observations with the results of  [L1] and \S 3 to draw some global
conclusions.

Given $\widehat\ii\subset S^2T^*$, a system of quadrics arising from the second
fundamental form at a general point of a variety,
  any cubic in the third
 fundamental form must be in  its prolongation,
$\widehat\ii\up 1 :=
(\widehat\ii\ot T^*)\cap S^3T^*$  (see [GH],(1.47), [L1],(3.12)). We will now
define a complement to
$\widehat\ii\up 1  $ in $\widehat\ii\ot T^*$  consisting of the space of
linear syzygies.
Recall that as a $Gl(T^*)$ module, $T^{*\ot 3}$ naturally splits into
three  factors; $S^3T^* \oplus \Lambda^3T^*
\oplus (S^{(21)}T^*)^{\oplus 2}$ where the last factor is two copies of
the irreducible $Gl(T^*)$ module obtained from the Young diagram with two boxes
in the first row and one in the second (hence the notation $(21)$).
We can choose two such copies as follows: Let
$$
\align
&S^{(21)}_{\circ}T^* := \tker (S^2T^*\ot T^*\ra S^3T^*)\tag 6.1\\
&S^{(21)}_{\Lambda}T^*:=\tker (\Lambda^2T^*\ot T^*\ra \Lambda^3T^*)
\tag 6.2
\endalign
$$
Then
$$
S^2T^*\ot T^* = S^3T^*\oplus \schur T^*\tag 6.3
$$
Given $A\subset S^2T^*$, define
$$
A^{[1]}:= (A\ot T^*)\cap\schur T^*\tag 6.4
$$
so
$$
A\ot T^* = A\up 1  \oplus A^{[1]}.\tag 6.5
$$

If $\widehat\ii\up 1   =0$, then $III_X=0$, and if
$\widehat\ii^{[1]}=0$ and $\widehat\ii\up 1   =0$,  then by (4.18)
the space of quadrics containing $X$ is at most $a$ dimensional.

 We now show that if the codimension of $X$ is
sufficiently small, then
$\widehat\ii^{[1]}=0$ and $\widehat\ii\up 1   =0$;
i.e. that
$\widehat\ii\ot T^*$ is disjoint form the two
$Gl(T)$ invariant linear spaces in $S^2T^*\ot T^*$.
\proclaim{Theorem 6.6}  Let  $X^n\subset\bpp{n+a}$ be a variety
and let $x\in X$ be any smooth point.   Let $b=\tdim X_{sing}$. (Set $b=-1$ if
$X$ is smooth.)  If
$a< \frac{n+1-(b+1)}{2}$, then $\widehat\ii_x\up 1  =0$, where
$\widehat\ii_x\up 1 =  S^3T^*\cap (\widehat\ii_x \ot T^*) $
is the {\it prolongation} of $\widehat\ii_x$.
\endproclaim
\demo{Proof}
Say there were a nonzero
$P\in\widehat\ii\up 1  $.  Consider
$$
\intprod P : T\ra \widehat\ii .\tag 6.7
$$
Applying ([L1],(6.1)) with  $A_x= \text{(image}\intprod P)$ gives
$$
\text{dim(Singloc(image}\intprod P))\leq
2(a-1) - (\text{dim(image}\intprod P)-1) + (b+1).\tag 6.8
$$
Observe that
$$
n=\text{dim(ker}(\intprod P ))+ \text{dim(image}(\intprod P))\tag 6.9
$$
and
$$
\text{ker}(\intprod P)\subseteq\text{ Singloc(image}(\intprod P ))\tag 6.10
$$
so
$$
n-\tdim (\text{image}(\intprod P))\leq 2a -1-\tdim (\text{image}(\intprod P))
+b+1\tag 6.11
$$
i.e.
$$
a\geq \frac{n+1-(b+1)}2\qed\tag 6.12
$$
\enddemo

\noindent{\bf Remark 6.13}. (6.6) gives a new proof of ([L1], (4.15)) stated in
the introduction
with refined information about what the structure of  the second and third
fundamental forms
of a variety must be in small codimension.

To study $\widehat\ii^{[1]}$ is a bit more difficult.
An $R\in \widehat\ii^{[1]}$ is a map $R: T\ra\widehat\ii $ and
$\ker R\subseteq$Baseloc$\widehat \ii$ so if we had a theorem directly
restricting the size of base loci of subsystems instead of singular
loci, we would be in better shape.

In general,  we can think of $II_X$ as
$$
\align ii: T&\ra T^*\ot N\tag 6.14\\
w&\mapsto II_X(w,\cdot )\endalign
$$
If we restrict $ii$ to $\tker R$, we get a map
$$
ii': \tker R \ra (\tker R)\upperp\ot N_R\tag 6.15
$$
where $N_R :=ii(\tker R) (T^*)$.
Now  $R\in \tker (\widehat\ii\ot T^*\ra S^3T^*)$
implies $(\tker R)\upperp \simeq N_R$ and moreover
that $ii'$ descends to a map
$$
ii' :\tker R \ra \Lambda^2 (\tker R)\upperp \simeq \Lambda^2 N_R.\tag 6.16
$$

In bases, the situation is as follows:
Let $L := \tker R\upperp \subset T^*$ and let $M\subset T^*$ be any complement
to $L$.
We may write $S^2T^* = S^2 L \oplus (L\ot M) \oplus S^2M$.

Let $A=R(T)\subset\widehat\ii$. Let $\{ q^j\}$ be a basis of $A$
and $\{ l^j \}$ a basis of $L$.
Write $q^j= b^j + c^j + d^j$ reflecting the decomposition of $S^2T^*$.
Then
$$
\align
 &b^j= b^i_{jk}l^jl^k, c^j=m^j_kl^k,\ d^j=0\tag 6.17\\
&\text{ with } b^i_{jk}\in\Bbb C,\ \frak S_{ijk}b^i_{jk}=0, \
m^i_j=-m^j_i\in M.\endalign
$$

By rechoosing $M$ if necessary, we may also assume
$$
b^i_{jk}=b^j_{ik} \tag 6.18
$$

\proclaim{Lemma 6.19}Let $A^p\subset S^2T^*$ be an $p$-dimensional
system of quadrics on an $n$ dimensional vector space.  Say there
is a linear syzygy
$$
l^1Q_1+\hdots + l^pQ_p = 0
$$
where both
$l^i\in T^*$ and $Q_i\in A$ are  independent sets of vectors.
Then $\forall Q\in A$,
$$
\text{rank }Q\leq 2(p-1).
$$
\endproclaim
\demo{Proof}
Let
$$
Q= \lambda^1Q_1 +\hdots + \lambda^pQ_p, \ \ \lambda^i\in\Bbb C.\tag 6.20
$$
Let $l_p\hd l_p$ be a dual basis to $l^1\hd l^p$ and
using (6.17), (6.18), observe
that
$$
(\lambda^1l_1+\hdots + \lambda^pl_p)\intprod Q=0.\tag 6.21
$$
On the other hand, we may write
$$
Q= l^1\alpha^1+\hdots l^p\alpha^p \tag 6.22
$$
for some $\alpha^1\hd \alpha^p\in T^*$. Now change bases
in $L$ such that
$l^p=(\lambda^1l_1+\hdots \lambda^pl_p)^*$. Then $\alpha^p=0$ and
we see $\trank Q\leq 2(p-1)$\qed\enddemo

Taking the $m_{ij}$ all independent (for $i<j$), produces quadrics of rank
$2(p-1)$, so the result is sharp.

Applying ([L1],(6.1)) to (6.19),(4.18) yields:

\proclaim{Lemma 6.23}Let  $X^n\subset\bpp{n+a}$ be a variety
and let $x\in X$ be a general point.   Let $b=\tdim X_{sing}$. (Set $b=-1$ if
$X$ is smooth.)
If
$a<\frac{n-(b+1)+3}{3}
$,
then there are no linear syzygies  in $\ii_x$.
\endproclaim
\demo{Proof} Lemma 6.20  says that if there is such a
syzygy, then there is a $p$ dimensional subsystem  $A_x\subseteq\ii_x$ with
the property that  no quadric in $A_x$ is of rank greater than $2(p-1)$.
Now  ([L1],(6.1)) applied to $A_x$ yields
$$
2(a-1) + (b+1) \geq n+p-1-2(p-1) \tag 6.24
$$
i.e.
$$
2(a-1) + p\geq n-(b+1)+1 \tag 6.25
$$
 Finally just observe that $p\leq a$.\qed
\enddemo

Applying (6.23) to (4.18) yields:

\proclaim{Theorem 6.26 }
Let  $X^n\subset\bpp{n+a}$ be a variety and $x\in X$ a general point.
  Let $b=\tdim X_{sing}$. (Set $b=-1$ if $X$ is smooth.)
If
$a<\frac{n-(b+1)+3}{3}
$
then
$$
\align &
\tdim \{ \text {quadrics osculating to order three at } x \}\leq
a
+\binom{a+1}2-1 \tag 6.27\\
&
\tdim \{ \text {quadrics osculating to order four at } x \}\leq
a-1.
\endalign
$$
 Equality occurs in the first (respectively second)  line of  (6.27)
 if and only if
(4.9.3) (resp. (4.9.4))   hold at $x$.
  If  the generalized Monge system
(4.17)  holds, then  $X$ is a complete intersection of the $a-1$
dimensional family of quadrics osculating to order four.
\endproclaim

 \demo{Proof}
It only remains to show there cannot be any
hypersurfaces of higher degree generating new elements of the
ideal of $X$.
The only way to have an equation of higher degree
is if there exists a nontrivial polynomial
$P(Q^1\hd Q^a)=0$, which is not possible in this codimension range as the
existence
of such a polynomial would imply the secant variety of $X$ is degenerate.\qed
\enddemo

\proclaim{Corollary 6.28}Let  $X^n\subset\bpp{n+a}$ be a variety and
 $x\in X$ a general point.
  Let $b=\tdim X_{sing}$. (Set $b=-1$ if $X$ is smooth.)
If
$a<\frac{n-(b+1)+3}{3}
$
then any quadric osculating to order four at $x$ is smooth at $x$ and
any quadric osculating to order five at $x$ contains $X$.
\endproclaim

\proclaim{Corollary 6.29 }
Let  $X^n\subset\bpp{n+a}$ be a variety
with $I_X$ generated by quadrics.   Let $b=\tdim
X_{sing}$. (Set $b=-1$ if $X$ is smooth.)
If
$a<\frac{n-(b+1)+3}{3},
$
then $X$ is a complete intersection.
\endproclaim

\bigpagebreak

{\bf References}:

\smallpagebreak

\noindent [BCG${}^3$] Bryant, R., Chern, S., Gardner, R., Goldschmidt,
Griffiths, P.,
 {\it Exterior Differential Systems}, Springer-Verlag, (1991) 475 pp.
{}.

\smallpagebreak

\noindent [C] Cartan, E., {\it  La methode du repere mobile},
Oeuvres Completes; Vol. 3, Part. 2, 1955,
pp. 1259-1320.

\smallpagebreak

\noindent [E] Ein, L., {\it  Varieties with
small dual varieties, I}, Invent. math {\bf
86}(1986) pp. 63-74.

\smallpagebreak

\noindent [F] Fubini,G., {\it Il problema della deformazione proiettiva
delle ipersuperficie}, Rend. Acad. Naz. dei Lincei {\bf 27}(1918) pp. 147-155.

\smallpagebreak

\noindent [FL] Fulton, W. and R. Lazarsfeld, {\it Connectivity and its
applications in algebraic geometry}, in Algebraic Geometry (Proceedings). Lect.
Notes in Math., No. {\bf 862}, pp. 26-92.

\smallpagebreak

\noindent [GH] Griffiths, P.A. and J. Harris, {\it Algebraic Geometry and Local
Differential Geometry}, Ann. scient. Ec. Norm. Sup.  {\bf 12} (1979), pp.
355-432.

\smallpagebreak

\noindent [H] Hartshorne, R., {\it Varieties of small codimension in projective
space}, Bull. A.M.S. {\bf 80} (1974), pp. 1017-1032.

\smallpagebreak

\noindent [JM] Jensen, G and E. Musso, {\it
Rigidity of hypersurfaces in complex projective space}, Ann. scient. Ec. Norm.
Sup.
 {\bf 27} (1994), pp. 227-248.

\smallpagebreak

\noindent [L1] Landsberg, J.M. {\it On Second Fundamental Forms
of Projective Varieties}.   Invent. Math.  {\bf 117}(1994) pp. 303-315.

\smallpagebreak

\noindent [L2] Landsberg, J.M. {\it On Degenerate Secant
and Tangential Varieties and Local Differential Geometry}.  Preprint available.
alg-geom  9412012

\smallpagebreak

\noindent [LVdV]  Lazarsfeld, R. and A. Van de Ven {\it Topics in the Geometry
of Projective Space, Recent Work of F.L. Zak}, DMV Seminar (1984), Birkhauser.

\smallpagebreak

\noindent [L'v] L'vovsky, S. {\it On Landsberg's criterion for
complete intersections} alg-geom 9408006.

\smallpagebreak

\noindent [S] Spivak, M. {\it A comprehensive introduction to
differential geometry}, Publish or Perish (1979), Vols. I-V.

\smallpagebreak

\noindent [T]  Terracini, {\it Alcune questioni sugli spazi tangenti e
osculatori ad una varieta, I, II, III}. Atti
Della Societa dei
Naturalisti e Matematici, (1913) pp. 214-247.

\smallpagebreak

\noindent [Z]  Zak, F.L., {\it Tangents and Secants of Algebraic Varieties},
AMS Translations of mathematical monographs, {\bf 127}(1993).

\enddocument